\documentclass[smallextended,numbook,runningheads]{svjour3}     
\smartqed  
\usepackage{graphicx}
\usepackage{amsmath}
\usepackage{epstopdf} 
\usepackage{mathptmx}      
\usepackage{amsmath}
\usepackage{mathtools, bm}
\usepackage{amssymb, bm}
\usepackage{cite}
\usepackage{mathrsfs}
\usepackage{lmodern}
\usepackage{multirow}
\usepackage{subfigure}
\usepackage{color}
\journalname{BIT}
\begin{document}

\title{Arbitrary high order A-stable and B-convergent numerical methods for  ODEs via deferred correction\thanks{The authors would like to acknowledge the financial support of the Discovery Grant Program of the Natural Sciences and Engineering Research Council of Canada (NSERC) and a scholarship to the first author from the NSERC CREATE program ``G\'enie par la Simulation''. 
}}

\titlerunning{Arbitrary order A-stable methods for ODE}        

\author{Saint-Cyr E.R.\ Koyaguerebo-Im\'e,\\ Yves Bourgault}


\institute{Saint-Cyr E.R.\ Koyaguerebo-Im\'e \and Yves Bourgault\thanks{Department of Mathematics and Statistics,
         University of Ottawa, STEM Complex, 150 Louis-Pasteur Pvt, Ottawa, ON, Canada, K1N 6N5, Tel.: +613-562-5800x2013 
  (\email{skoya005@uottawa.ca}, \email{ybourg@uottawa.ca}
}}

\date{Received: date / Accepted: date}

\maketitle

\begin{abstract}
 This paper presents a sequence of deferred correction (DC) schemes built recursively from the implicit midpoint scheme for the numerical solution of general first order ordinary differential equations (ODEs). It is proven that each scheme is A-stable, satisfies a B-convergence property, and that the correction on a scheme DC2j of order 2j of accuracy leads to a scheme DC2j+2 of order 2j+2. The order of accuracy is guaranteed by a deferred correction condition. Numerical experiments with standard stiff and non-stiff ODEs are performed with the DC2, ..., DC10 schemes. The results show a high accuracy of the method. The theoretical orders of accuracy are achieved together with a satisfactory stability.
\end{abstract}

\section{Introduction}
\label{intro}
In \cite{MR2058857,kress2002deferred}, Gustafsson and  Kress introduced a new version of deferred correction (DC) strategy for the numerical solution of linear systems of ordinary differential equations (ODE) \cite{MR2058857} and  initial boundary value problems \cite{kress2002deferred}, under  a monotonicity condition. Numerical experiments with one-dimensional linear parabolic and hyperbolic equations were performed and showed that the method is effective (orders 2, 4 and 6 of accuracy are achieved). We propose to extend the method from \cite{MR2058857,kress2002deferred} to the time-discretization of more general time-evolution partial differential equations (PDEs). In this paper, we restrict to the case of the initial value problem (IVP)
\begin{equation}
\label{a24}
\left\lbrace 
\begin{array}{cccc}
\displaystyle \frac{d u}{dt}&=&F(t,u),&~~ t\in [0,T],\\
u(0)&=&u_0,&
\end{array}
\right.
\end{equation}
where the unknown $u$ is from $[0,T]$ into a Banach space $X$, $u_0$ is a given data and $F$ is a sufficiently differentiable function such that $u$ exists and is  sufficiently differentiable. The main objective is to show the properties of the numerical method (consistency, stability, convergence and order of accuracy). A complete analysis of the DC method applied to reaction-diffusion equations leads to an arbitrary high order and unconditionally stable method (see \cite{koyaguerebo2020unconditionally}). 

The DC method is used to improve the order of accuracy of numerical methods of lower order. This method is explored by many authors, e.g. \cite{schild1990gaussian,auzinger2016encyclopedia,MR2058857,kushnir2012highly,hansen2011order,dutt2000spectral,daniel1967interated,IntegralDC2010}. The method in  \cite{daniel1967interated} is an application of iterative deferred correction (IDC). The authors proved that an asymptotic improvement of order $p$ can be accomplished, from a scheme of order $p$, at  each step of the IDC procedure, provided suitable finite difference operators are employed. Numerical experiments are performed with the IDC applied to the trapezoidal rule, Taylor-2 and Adams-Bashforth of order 2. The results are promising even though they point out some difficulties of the proposed algorithms: inaccuracy for ``large'' time step and no asymptotic improvement for high levels of correction. The approaches in \cite{kushnir2012highly,hansen2011order,dutt2000spectral,auzinger2016encyclopedia,MR2058857,IntegralDC2010} are quite similar and consist in a linear perturbation of a low order scheme. However, solving stiff problems (problems extremely hard to solve by standard explicit methods \cite{spijker1996stiffness}) is a challenge unfavorable for these methods. In particular, the method in \cite{kushnir2012highly}, concerning a highly accurate solver for stiff ODEs, requires sufficiently small time
steps for moderately stiff problems while convergence is reduced to order 2 for
``very stiff'' problems.

Our schemes are based on nonlinear perturbations (corrections) of the implicit midpoint rule and inherit the A-stable property of the trapezoidal rule \cite{MR0170477} at any stage of the correction. Starting from an approximation $\left\lbrace u^{2,n}\right\rbrace _{n=0}^N$ of the exact solution $u$ by the implicit midpoint rule on a uniform partition $0=t_0<t_1<\cdots< t_N=T$ of $[0,T]$, at the stage $j=1,2,\cdots$ of the correction we obtain an approximation $\left\lbrace u^{2j+2,n}\right\rbrace _{n=0}^N$ of $u$, expected to be of order $2j+2$ of accuracy, on the same partition. Each approximate solution $\left\lbrace u^{2j,n}\right\rbrace _{n=0}^N$ to be corrected is subject to a deferred correction condition (DCC) which guarantees the improvement of the order of accuracy. We prove that if $\left\lbrace u^{2j,n}\right\rbrace _{n=0}^N$  satisfies the DCC and its correction $\left\lbrace u^{2j+2,n}\right\rbrace _{n=0}^N$ converges to $u$ at the discrete points $0=t_0<t_1<\cdots< t_N=T$ (or is simply bounded, when $X$ is finite dimensional) then $\left\lbrace u^{2j+2,n}\right\rbrace _{n=0}^N$  approximates $u$ with order $2j+2$. Moreover, provided the function $F$ is Lipschitz with respect to its second variable or satisfies a one-sided Lipschitz condition, each $\left\lbrace u^{2j,n}\right\rbrace _{n=0}^N$ satisfies the DCC and then converges with order $2j$ of accuracy, for arbitrary positive integer $j$. We also prove that each DC scheme involving $\left\lbrace u^{2j,n}\right\rbrace _{n=0}^N$ is $B$-stable. The theory is illustrated by numerical tests, for the schemes of order 2, 4, ..., 10.

The paper is organized as follows: in section 2 we recall some basic results from finite difference approximations and present the DC schemes; section 3 deals with the consistency of the method; the analysis of convergence and order of accuracy together with a B-convergence result are given in section 4; absolute stability is proved is section 5, and section 6 is devoted to numerical experiments.

\section{Deferred correction schemes for the implicit midpoint rule}
\label{sec:main}
We suppose that $\displaystyle F\in C^{2p+2}\left( [0,T]\times X,X \right)$, for a positive integer $p$, so that (\ref{a24}) has a unique solution $\displaystyle u\in C^{2p+3}\left( [0,T],X \right)$. We simply denote by $\|\cdot\|$, the norm in the Banach space $X$. For a time step $k>0$, we denote $t_n=nk$ and $t_{n+1/2}=(n+1/2)k$, for each integer $n$. This implies that $t_0=0$. We consider the time steps $k$ such that $0=t_0<t_1<\cdots< t_N=T$ is a partition of $[0,T]$, for a non-negative integer $N$. The centered, forward and backward difference operators $D$, $D_+$ and $D_-$, respectively, related to $k$ and applied to $u$, are defined as follows:
$$
Du(t_{n+1/2})=\frac{u(t_{n+1})-u(t_n)}{k},$$
$$D_+u(t_{n})=\frac{u(t_{n+1})-u(t_n)}{k},$$ 
and 
$$D_-u(t_{n})=\frac{u(t_{n})-u(t_{n-1})}{k}, n\geq 1.$$ The average operator is denoted by $E$:
$$E u(t_{n+1/2})=\widehat{u}(t_{n+1})=\frac{u(t_{n+1})+u(t_n)}{2}.$$
The composition of $D_+$ and $D_-$ is defined recursively. They commute, that is $(D_+D_-)u(t_n)=(D_-D_+)u(t_n)=D_-D_+u(t_n)$, and satisfy the identities
\begin{equation}
\label{bb1}
(D_+D_-)^mu(t_n)=k^{-2m} \sum_{i=0}^{2m}(-1)^i {{2m}\choose {i}}u(t_{n+m-i} ),
\end{equation}and
\begin{equation}
\label{bb2}
D_-(D_+D_-)^mu(t_n)=k^{-2m-1}\sum_{i=0}^{2m+1}(-1)^i{{2m+1}\choose {i}}u(t_{n+m-i}),
\end{equation}
for each integer $m\geq 1$ such that $0\leq t_{n-m-1}\leq t_{n+m}\leq T$. We have the estimate
\begin{equation}
\label{b3}\left \Vert D_+^{m_1}D_-^{m_2}u(t_n) \right \Vert\leq \max_{0\leq t\leq  T}\left \Vert\frac{d^{m_1+m_2}u}{dt^{m_1+m_2}}(t)\right \Vert,
\end{equation}
provided $[t_{n-m_2},t_{n+m_1}] \subset [0,T]$ and $m_1+m_2\leq 2p+3$ (see \cite[p.249]{isaacson1966analysis} or \cite{koyaguerebo2020finite}).

If $\left\lbrace u^n \right\rbrace_n$ is a sequence of approximation of $u$ at the discrete points $t_n$, the finite difference operators apply to $\left\lbrace u^n \right\rbrace _n$, and we define
$$
Du^{n+1/2}=D_+u^{n}=D_-u^{n+1}=\frac{u^{n+1}-u^n}{k},$$and
$$E u^{n+1/2}=\widehat{u}^{n+1}=\frac{u^{n+1}+u^n}{2}.$$

From the centered finite difference approximation (see \cite[Thm 5]{koyaguerebo2020finite} or \cite{hildebrand1974introduction,chung2010computational,dahlquist2008numerical}) we have 
\begin{equation}
\label{bb3}
\frac{d u}{dt}(t_{n+1/2})=
\frac{u(t_{n+1})-u(t_n)}{k}-\sum_{i=1}^jc_{2i+1}k^{2i}(D_+D_-)^iDu(t_{n+1/2}))+O(k^{2j+2}) 
\end{equation}
and
\begin{equation}
\label{bb4}
u(t_{n+1/2}) =\frac{u(t_{n+1})+u(t_n)}{2}-\sum_{i=1}^jc_{2i}k^{2i}(D_+D_-)^iE u(t_{n+1/2})+O(k^{2j+2}),
\end{equation} 
for each integer $j=1,2,\cdots,p$. These approximations lead to the schemes 
\begin{equation}
\label{a25}
\begin{aligned}
\frac{u^{n+1}-u^n}{k}-&\sum_{i=1}^jc_{2i+1}k^{2i}(D_+D_-)^iDu^{n+1/2}\\&=F\left( t_{n+1/2},\frac{u^{n+1}+u^n}{2}-\sum_{i=1}^jc_{2i}k^{2i}(D_+D_-)^iE u^{n+1/2}\right).
\end{aligned}
\end{equation}
The schemes (\ref{a25}) are multi-steps and prone to stability restrictions. We resort to DC method to transform them into a sequence of one step schemes as follows: 
{For $j=0$, we have the implicit midpoint rule}
\begin{equation}
\label{a26}\frac{u^{2,n+1}-u^{2,n}}{k}=F\left( t_{n+1/2},\frac{u^{2,n+1}+u^{2,n}}{2}\right),~~ u^{2,0}=u_0.
\end{equation}

{For $j\geq 1$,}
\begin{equation}
\begin{aligned}
\label{a27}
 &\frac{u^{2j+2,n+1}-u^{2j+2,n}}{k}-\sum_{i=1}^jc_{2i+1}k^{2i}(D_+D_-)^iDu^{2j,n+1/2}\\&=F\left( t_{n+1/2},\frac{u^{2j+2,n+1}+u^{2j+2,n}}{2}-\sum_{i=1}^jc_{2i}k^{2i}(D_+D_-)^iEu^{2j,n+1/2}\right),
\end{aligned}
\end{equation}
\begin{equation}
\label{a28}u^{2j+2,0}=u_0.
\end{equation}
The scheme (\ref{a27})-(\ref{a28}) has unknowns $u^{2j+2,n}$, $n=1,2,...,N$, and is deduced from (\ref{a25}) by substituting the unknown $u^n$ under the summation symbols by $u^{2j,n}$. The index $2j$ indicates that $\left\lbrace u^{2j,n}\right\rbrace _n$ is expected to be an approximation of the exact solution $u$ with order $2j$ of accuracy.  We call the schemes (\ref{a27})-(\ref{a28}) Deferred Correction of order $2j+2$ for the implicit midpoint rule, denoted DC2j+2.

\begin{remark}
\label{rmk:b1}
The scheme (\ref{a27})-(\ref{a28}), for $n=1,2,3,\cdots,j$, should involve unknowns $u^{2j,-1}, ... , u^{2j,-j}$ which represent approximate solutions of (\ref{a24}) at the discrete points $t=-k,...,-jk$, respectively. To avoid those approximations for $t<0$, we propose the following scheme which is efficient for the computation of $u^{2j+2,1},..., u^{2j+2,j}$, using only  points within the solution interval $[0,T]$.
\begin{equation}
\label{a73}
\begin{aligned} &\frac{u^{2j+2,n+1}-u^{2j+2,n}}{k}-k^{-1}\sum_{i=1}^j c_{2i+1}^jk_j^{2i+1}(D_+D_-)^iD \bar u^{2j,(2j+1)n+j+1/2}\\&=F\left( t_{n+1/2},Eu^{2j+2,n+1/2}-\sum_{i=1}^j c_{2i}^jk_j^{2i}(D_+D_-)^iE\bar u^{2j,(2j+1)n+j+1/2}\right),
\end{aligned}
\end{equation}
\begin{equation}
\label{a74b}u^{2j+2,0}=u_0.
\end{equation}
The finite difference operator in (\ref{a73}) are related to the time step $k_j=k/(2j+1)$. The approximations $\left\lbrace \overline{u}^{2j,m}\right\rbrace _m$ and $\left\lbrace u^{2j,n}\right\rbrace _n$ are computed from the same scheme, (\ref{a26}) or (\ref{a27})-(\ref{a28}), but for the time steps $k_j$ and $k$, respectively. The scheme (\ref{a73}) results from the finite difference approximations
\begin{equation}
\label{01} 
u'(t_{n+1/2})=
\frac{u(t_{n+1})-u(t_n)}{k}-\frac{1}{k}\sum_{i=1}^jc^j_{2i+1}k_j^{2i+1}D(D_+D_-)^iu(\tau_{j+1/2}) +O(k_j^{2j+2})
\end{equation}
and
\begin{equation}
\label{02} 
u(t_{n+1/2}) =\frac{u(t_{n+1})+u(t_n)}{2}-\sum_{i=1}^jc^j_{2i}k_j^{2i}(D_+D_-)^iE u(\tau_{j+1/2})+O(k_j^{2j+2}),
\end{equation} where $t_n=\tau_0<\tau_1<...<\tau_{2j+1}=t_{n+1}$, with $\tau_m=t_n+mk_j$, for $m=1, 2,\cdots , 2j+1$.
Table \ref{tab:3} gives the coefficients ${c}^j_i$ for $j=1,2,3,4$.

\begin{table}[!ht] 
\caption{Coefficients of the approximations (\ref{01})-(\ref{02}) for $j=1,2,3,4$}
\label{tab:3}
 \centering      
\begin{tabular}{lllllllllll}  
\hline\noalign{\smallskip} 
   $j$& ${c}^j_2$ &~${c}^j_3$ &~~~ ${c}^j_4$ &~~~~${c}^j_5$& ~${c}^j_6$&~~${c}^j_7$&~~~~${c}^j_8$&~~~~~${c}^j_9$&~~~\\[0.75ex]
  1&  ~$\frac{9}{8}$        &~$\frac{9}{8}$       \\[0.75ex]
  2&  $\frac{25}{8}$        &$\frac{125}{24}$         &$\frac{125}{128}$      &$\frac{125}{128}$     \\[0.75ex]
  3&  $\frac{49}{8}$        &$\frac{343}{24}$         &$\frac{637}{128}$      &$\frac{13377}{1920}$     &$\frac{1029}{1024}$      &$\frac{1029}{1024}$      \\[0.75ex]
  4&  $\frac{81}{8}$        &$\frac{243}{8}$         &$\frac{1917}{128}$      &$\frac{17253}{640}$     &$\frac{7173}{1024}$      &$\frac{64557}{7168}$      &$\frac{32733}{32768}$   &$\frac{32733}{32768}$  \\[.75ex]
\noalign{\smallskip}\hline
\end{tabular} 
\end{table}  
\end{remark}

\begin{remark}
\label{rmk:b2} Each $u^{2j+2,n+1}$, $n\geq j$, is an iterative solution of the system
\begin{equation}
\label{o4} x-a_n^j-kF(t_{n+1/2},0.5x+b_n^j)=0,
\end{equation}
where $x$ is the unknown, and $a_n^j$ and $b_n^j$ are constants depending on $u^{2j+2,n}$ and $u^{2j,n+1+j},u^{2j,n+j},\cdots,u^{2j,n-j}$. The total number of vectors (in the solution space $X$) stored for the computation of $u^{2j+2,n+1}$ is $j^2+3j+1$: $u^{2j+2,n}$ and the $u^{2i,q}$, for $i=1,2,\cdots,j$, and $n+(j-i+1)(j+i)/2-2i \leq q\leq n+1+(j-i+1)(j+i)/2$.
\end{remark}

\begin{remark} 
From Remark \ref{rmk:b2}, only the implicit midpoint rule, DC2, is an implicit Runge-Kutta (RK) methods. Starting with DC4, all the DC2j methods of the form (\ref{a27})-(\ref{a28}) are not RK methods. For instance, $u^{4,n+1}$ depends on $u^{4,n}$ and some of the $u^{2,i}$, which $u^{2,i}$ evolve independently and are not stages computed from $u^{4,n}$. As we will see in Section \ref{sec:Astability}, the analysis of A-stability, in particular the proof of lemma \ref{lem:4}, shows that it is impossible to write a recurrence $u^{2j+2,n+1}=R(z)\,u^{2j+2,n}$ from (\ref{a27}) when $j\ge 1$, as one would get by applying any RK method to Dahlquist equation. This is the main ingredient behind the A-stability of our DC2j methods independently of the order of accuracy.
\end{remark}

\section{Deferred correction condition (DCC)}
\label{sec:DCC}
In this section we give a sufficient condition for the scheme (\ref{a27})-(\ref{a28}) to achieve order $2j+2$ of accuracy. Hereafter, the letter $C$ will denote any constant independent from $k$, and that can be calculated explicitly in terms of known quantities. The exact value of $C$ may change. We have the following definition:
 
\begin{definition} {\bfseries (Deferred Correction Condition)}
\label{def:1} Let $u$ be the exact solution of the Cauchy problem  (\ref{a24}). Given a positive integer $j$, a sequence $\left\lbrace u^{2j,n}\right\rbrace _{n=0}^N$ of approximations of $u$, at the discrete points $0=t_0<\cdots <t_N=T$,  is said to satisfy the Deferred Correction Condition $(DCC)$ for the implicit midpoint rule if $\left\lbrace u^{2j,n}\right\rbrace _{n=0}^N$ approximates $u$ with  order $2j$ of accuracy, and we have
\begin{equation}
\label{a30} 
 \Vert (D_+D_-)D(u^{2j,n+1/2}-u(t_{n+1/2}))\Vert+\Vert D_+D_-(u^{2j,n+1}-u(t_{n+1}))\Vert \leq C k^{2j},
\end{equation}for $n=1,2,...,N-2$ and $k\leq k_0$, where $k_0>0$ is fixed and $C$ is a constant independent from $k$.
\end{definition}

\begin{remark}
\label{rmk:4} Condition (\ref{a30}) is equivalent to
\begin{equation}
\label{cc1}
\left \| \sum_{i=1}^jc_{2i}k^{2i}\left( D_+D_-\right)^i\left(u^{2j,n}-u(t_n) \right) \right \|\leq Ck^{2j+2},
\end{equation}
and
\begin{equation}
\label{cc2}
\left \| \sum_{i=1}^j(c_{2i+1}-c_{2i})k^{2i}\left( D_+D_-\right)^iD\left(u^{2j,n+1/2}-u(t_{n+1/2}) \right) \right \|\leq Ck^{2j+2},
\end{equation}
for $n=j,j+1,\cdots , N-j$. This is due to the transform
$$k^{2i}\left(D_+D_-\right)^i\left(u^{2j,n}-u(t_n) \right)=k^2\sum_{l=0}^{i-1}(-1)^l{{2i-2}\choose{l}}D_+D_-\left(u^{2j,n}-u(t_n) \right) $$
and a similar transform for $k^i\left(D_+D_-\right)^iD\left(u^{2j,n+1/2}-u(t_{n+1/2}) \right)$.
\end{remark}

We have the following result:

\begin{theorem}
\label{theo:3} Let $u$ be the exact solution of (\ref{a24}) and $\left\lbrace u^{2j,n}\right\rbrace _{n=0}^N$, $j=1,\dots,p$, a sequence of approximations of $u$ satisfying DCC for the implicit midpoint rule. Let $\left\lbrace u^{2j+2,n}\right\rbrace _{n=0}^N$ be the solution of (\ref{a27})-(\ref{a28}) built from $\left\lbrace u^{2j,n}\right\rbrace _{n=0}^N$. We suppose that $u^{2j+2,1},...,u^{2j+2,j}$ are given and satisfy 
\begin{equation}
\label{a29}
\Vert u^{2j+2,n}-u(t_n) \Vert\leq Ck^{2j+2},~~\mbox{ for }~n=1,2,...,j,
\end{equation}
 where $C$ is a constant independent from $k$. Furthermore, we suppose that one of the following four conditions holds:

\item[(i)] $F$ is Lipschitz with respect to the second variable $x$: there exists $\mu \geq 0$ such that
\begin{equation}
\label{a31a} \|F(t,x)-F(t,y)\|\leq \mu\|x-y\|,~~\forall (t,x,y)\in [0,T]\times X\times X.
\end{equation}
\item[(ii)] $X$ is finite dimensional, and $\left\lbrace u^{2j+2,n}\right\rbrace _{n=0}^N$ remains close to $u$ in the sense that there exists  $M>0$ such that
\begin{equation}
\label{a32}
\Vert u^{2j+2,n}-u(t_n) \Vert \leq M,~~\mbox{ for each }~n=0,1,...,N.
\end{equation}
\item[(iii)]  $X$ is infinite dimensional, and $\left\lbrace u^{2j+2,n}\right\rbrace _n$ converges to the exact solution $u$.
\item[(iv)] $X$ is a Hilbert space with inner product $\left( .,. \right) $, and $F$ satisfies the following so-called one-sided Lipschitz condition, with a one-sided Lipschitz constant $\beta \in \mathbb{R}$:
\begin{equation}
\label{a31} \left( F(t,x)-F(t,y),x-y\right)\leq \beta\|x-y\|^2, ~~\forall (t,x,y)\in [0,T]\times X\times X. 
\end{equation}
Then $\left\lbrace u^{2j+2,n}\right\rbrace _n$ approximates  $u$ with order $2j+2$ of accuracy, that is 
\begin{equation}
\label{a32c}\|u^{2j+2,n}-u(t_n)\|\leq Ck^{2j+2},~~\mbox{ for each }~n=0,1,...,N,
\end{equation}
where $C$ is a constant depending only on $j$, $T$, DCC, a Lipschitz constant on $F$ and the derivatives of $u$ up to order $2j+3$, for time steps $k$ sufficiently small.
\end{theorem}

\begin{proof}
\item[1.]First we consider the case where the function $F=F(t,x)$ is Lipschitz with respect to the second variable $x$. Combining (\ref{a24}) and (\ref{a27}), we obtain the identity
\begin{equation}
\label{a33}
\begin{aligned}
&D\Theta^{2j+2,n+1/2} =\sigma^{2j+2,n+1/2}+ (\Lambda^j -\Gamma^j)D\left( u^{2j,n+1/2}-u(t_{n+1/2})\right) \\&+F\left( t_{n+1/2}, \widehat{u}^{2j+2,n+1 }-\Gamma^j \widehat{u}^{2j,n+1}\right)-F\left( t_{n+1/2},  \widehat{u}(t_{n+1})-\Gamma^j \widehat{u}(t_{n+1})\right)  ,
\end{aligned}
\end{equation}where $\Lambda^j$ and $\Gamma^j$ are finite difference operators defined for arbitrary integer $j \geq 1$ by
$$\Lambda^ju(t_n) =\sum_{i=1}^jc_{2i+1}k^{2i}(D_+D_-)^i u(t_n),$$
and
$$\Gamma^j u(t_n)=\sum_{i=1}^jc_{2i}k^{2i}(D_+D_-)^i u(t_n),$$
provided $u(t_{n\pm i})$ exists for $i=0,1,2,\cdots,j$. We have defined
\begin{equation}
\label{a34}
\Theta^{2j+2,n}=\left( u^{2j+2,n}-u(t_{n})\right) -\Gamma^j \left( u^{2j,n}-u(t_n)\right),
\end{equation}
and
\begin{equation*}
\begin{aligned}
\sigma^{2j+2,n+1/2}&=\left[ u'(t_{n+1/2})-Du(t_{n+1/2})+\Lambda^j Du(t_{n+1/2})\right]\\&-\left[ F(t_{n+1/2},u(t_{n+1/2}))-F( t_{n+1/2},  \widehat{u}(t_{n+1})-\Gamma^j \widehat{u}(t_{n+1}))\right] .
\end{aligned}
\end{equation*}
From (\ref{bb3}) we have
$$\left \Vert u'(t_{n+1/2})-Du(t_{n+1/2})+\Lambda^jD u(t_{n+1/2}) \right \Vert \leq Ck^{2j+2},$$
and, since $F$ is differentiable and $u$ is sufficiently regular, we deduce from the mean value theorem and the approximation (\ref{bb4}) that
\begin{equation*}
\left \Vert F(t_{n+1/2},u(t_{n+1/2}))-F( t_{n+1/2},  \widehat{u}(t_{n+1})-\Gamma^j \widehat{u}(t_{n+1})) \right \Vert \leq Ck^{2j+2}, 
\end{equation*}
for each $n=0,1,\cdots,N$, where $C$ is a constant depending only on $j$, $T$, a Lipschitz constant from $F$ and the derivatives of $u$ up to order $2j+3$. The last two inequalities imply that
\begin{equation}
\label{a36} 
\left \Vert\sigma^{2j+2,n+1/2}\right \Vert \leq Ck^{2j+2}.
\end{equation} Since the sequence $\left\lbrace u^{2j,n}\right\rbrace _n$ satisfies DCC, from Remark \ref{rmk:4} we have
\begin{equation}
\label{a37} 
\left \Vert \left( \Lambda^j -\Gamma^j\right) D\left( u^{2j,n+1/2}-u(t_{n+1/2})\right)  \right \Vert \leq Ck^{2j+2}.
\end{equation}
From the Lipschitz condition on $F$ we have
\begin{equation}
\label{a38} 
\begin{aligned}
&\left \Vert  F\left( t_{n+1/2}, \widehat{u}^{2j+2,n+1 }-\Gamma^j \widehat{u}^{2j,n+1}\right)-F\left( t_{n+1/2},  \widehat{u}(t_{n+1})-\Gamma^j \widehat{u}(t_{n+1})\right) \right  \Vert\\& \leq \mu\Vert \widehat{\Theta}^{2j+2,n+1}\Vert.
\end{aligned}
\end{equation}
Substituting inequalities (\ref{a36})-(\ref{a38}) in the identity (\ref{a33}), we deduce that 

$$\Vert D\Theta^{2j+2,n+1/2} \Vert \leq Ck^{2j+2}+\mu \Vert \widehat{\Theta}^{2j+2,n+1}\Vert, $$and it follows from the triangle inequality that

$$\Vert \Theta^{2j+2,n+1} \Vert \leq C\frac{k^{2j+3}}{2-\mu k}+\frac{2+\mu k}{2-\mu k}\Vert \Theta^{2j+2,n} \Vert ,$$ for $0\leq \mu k <2$. We then deduce by induction on $n$ that 

\begin{equation}
\label{a39a}
\Vert \Theta^{2j+2,n} \Vert \leq C\frac{1}{2-\mu k}\left( \frac{2+\mu k}{2-\mu k}\right)^{n-j-1} k^{2j+2}+\left( \frac{2+\mu k}{2-\mu k}\right)^{n-j}\Vert \Theta^{2j+2,j} \Vert .
\end{equation}
From hypothesis (\ref{a29}) and the DCC we have 
\begin{equation}
\label{a39}
\Vert \Theta^{2j+2,j} \Vert \leq \Vert u^{2j+2,j}-u(t_j)\Vert+ \left \Vert \Gamma^j(u^{2j,j}-u(t_j))\right  \Vert\leq Ck^{2j+2},
\end{equation}
where $C$ is a constant independent from $k$. Moreover, the sequence $\left\lbrace \left( \frac{2+\mu k}{2-\mu k}\right)^{n}  \right\rbrace_n $  is bounded above by  $\exp(2\mu T/(2-\varepsilon))$, for $0\leq \mu k \leq \varepsilon <2$. Whence 
$$\Vert \Theta^{2j+2,n} \Vert \leq Ck^{2j+2}.$$
Finally, by the triangle inequality, identity (\ref{a34}) and DCC, we have 

$$\Vert u^{2j+2,n}-u(t_n)\Vert \leq \Vert \Theta^{2j+2,n} \Vert+\left \Vert \Gamma^j (u^{2j,n}-u(t_n))\right \Vert \leq Ck^{2j+2},$$ where $C$ is a constant depending only on $j$, $T$, the DCC constant, $\mu$ and the derivatives of $u$ up to order $2j+3$. 
\newline 

\item[2.] Suppose that $\left\lbrace u^{2j+2,n} \right\rbrace_{n=0}^N $ satisfies (\ref{a32}) and $X$ is finite dimensional. We can write
\begin{equation*}
\begin{aligned}
F&\left( t_{n+1/2}, \widehat{u}^{2j+2,n+1 }-\Gamma^j \widehat{u}^{2j,n+1}\right)-F\left( t_{n+1/2},  \widehat{u}(t_{n+1})-\Gamma^j \widehat{u}(t_{n+1})\right)\\&=
 \int_0^1d_uF\left(t_{n+1/2},\widehat{u}(t_{n+1})-\Gamma^j \widehat{u}(t_{n+1})+s \widehat{\Theta}^{2j+2,n+1} \right)\left(\widehat{\Theta}^{2j+2,n+1}  \right) ds. 
 \end{aligned}
\end{equation*} 
From (\ref{a32}) and the DCC there exists $k_1>0$ such that $0<k\leq k_1\leq k_0$ implies 
$$\| \widehat{\Theta}^{2j+2,n+1} \|\leq M+Ck^{2j+2}\leq M+1. $$
On the other hand, we have

\begin{equation}
\label{saint1a}\|\widehat{u}(t_{n+1})-\Gamma^j \widehat{u}(t_{n+1})\|=\left \|\widehat{u}(t_{n+1})-\sum_{i=1}^j\sum_{l=0}^{2i}(-1)^{l}c_{2i}{{2i}\choose {l}}u(t_{n+i-l} )\right \|\leq R_{j+1},
\end{equation}
where
\begin{equation}
\label{saint1b}
R_{j+1}: =(j+1)\max_{0\leq t\leq T}\Vert u(t)\Vert \geq \left(1+\sum_{i=1}^j2^{2i}|c_{2i}|  \right)\max_{0\leq t\leq T}\Vert u(t)\Vert.
\end{equation}
It follows (\ref{a38}) for 
$$ \mu =\sup_{0\leq t\leq T,\Vert x\Vert \leq M+R_{j+1}+1}\left \Vert d_xF(t,x) \right \Vert.$$
Since $F$ is differentiable and the set $\left\lbrace x\in X: \Vert x\Vert \leq M+R_{j+1}+1 \right\rbrace $ is compact in the finite dimensional linear space $X$, the supremum exists and is finite.
The theorem is then deduced from the case (i).
\newline
\item[3.] If $\left\lbrace u^{2j+2,n}\right\rbrace _n$ converges to the exact solution $u$, 
taking the DDC and the finite difference formula (\ref{bb4}) into account, we have 
$$\left( \widehat{u}(t_{n+1})-\Gamma^j \widehat{u}(t_{n+1})+s \widehat{\Theta}^{2j+2,n+1}\right)  -u(t_{n+1/2})\rightarrow 0, \mbox{ as }k\rightarrow 0, \mbox{ for }0\leq s\leq 1.$$
It follows from the continuity of $u \mapsto d_uF(t,u)$ that there exists $0<k_2\leq k_0$ such that $0< k\leq k_2$  implies 
$$
\begin{aligned}
\Vert d_uF(t_{n+1/2},&\widehat{u}(t_{n+1})-\Gamma \widehat{u}(t_{n+1})+\tau \widehat{\Theta}^{2j+2,n+1} ) \Vert \leq  1+\max_{0\leq t\leq T}\Vert d_uF\left( t,u(t)\right)\Vert .
\end{aligned}$$
The theorem, in this case, follows by taking $\mu= 1+\max_{0\leq t\leq T}\Vert d_uF\left( t,u(t)\right)\Vert $ in (i).
\newline

\item[4.] Here we consider the case where $X$ is a Hilbert space and $F$ satisfies the monotonicity condition (\ref{a31}). Then, taking the inner product of the identity (\ref{a33}) with $\widehat{\Theta}^{2j+2,n+1}$, we deduce the inequality
\begin{equation}
\label{a40}
\begin{split}
\left( D\Theta^{2j+2,n+1/2},\widehat{\Theta}^{2j+2,n+1} \right) \leq 
\left( \sigma^{2j+2,n+1/2} ,\widehat{\Theta}^{2j+2,n+1}\right)+\beta \|\widehat{\Theta}^{2j+2,n+1}\|^2  \\\left( (\Lambda^j -\Gamma^j)D(u^{2j,n+1/2}-u(t_{n+1/2})),\widehat{\Theta}^{2j+2,n+1}\right) 
\end{split}
\end{equation}
since, according to (\ref{a31}), we have 
\begin{equation*}
\begin{aligned}
\left(  F\left(  t_{n+1/2}, \widehat{u}^{2j+2,n+1 }-\Gamma \widehat{u}^{2j,n+1}\right) -F\left(  t_{n+1/2},  \widehat{u}(t_{n+1})-\Gamma \widehat{u}(t_{n+1})\right) ,\widehat{\Theta}^{2j+2,n+1}\right) \\
\leq \beta \left \|\widehat{\Theta}^{2j+2,n+1}\right \|^2.
\end{aligned}
\end{equation*}
Inequalities (\ref{a36})-(\ref{a37}) together with the Cauchy-Schwartz inequality yield 
\begin{equation*}
\left \vert \left(\sigma^{2j+2,n+1/2},\widehat{\Theta}^{2j+2,n+1}\right)\right \vert \leq Ck^{2j+2}\|\widehat{\Theta}^{2j+2,n+1}\|,
\end{equation*}
and
\begin{equation*}
\left \vert \left( (\Lambda^j -\Gamma^j)D(u^{2j,n+1/2}-u(t_{n+1/2})),\widehat{\Theta}^{2j+2,n+1}\right)\right \vert \leq Ck^{2j+2}\|\widehat{\Theta}^{2j+2,n+1}\|,
\end{equation*}
where $C$ is a constant depending only on $j$, $T$, a Lipschitz constant on $F$ and the derivatives of $u$ up to order $2j+3$. Substituting the last three inequalities into (\ref{a40}), we obtain 
\begin{equation*}
\left( D\Theta^{2j+2,n+1/2},\widehat{\Theta}^{2j+2,n+1} \right) \leq Ck^{2j+2}\|\widehat{\Theta}^{2j+2,n+1}\|+\beta\|\widehat{\Theta}^{2j+2,n+1}\|^2 ,
\end{equation*} and we deduce from the identity
\begin{equation*}
 \left( D\Theta^{2j+2,n+1/2},\widehat{\Theta}^{2j+2,n+1} \right)=\frac{1}{2k}\left(\|\Theta^{2j+2,n+1}\|^2-\|\Theta^{2j+2,n}\|^2 \right)
\end{equation*}and the inequality
$$\|\widehat{\Theta}^{2j+2,n+1}\| \leq \frac{1}{2}\left(\|\Theta^{2j+2,n+1}\|+\|\Theta^{2j+2,n}\| \right)$$that 

$$\|\Theta^{2j+2,n+1}\|\leq C\frac{k^{2j+3}}{2-\beta k}+\frac{2+\beta k}{2-\beta k}\Vert \Theta^{2j+2,n} \Vert.$$
 The conclusion follows from the case (i), for $-2\leq \beta k<2$.
\end{proof}
\vspace{.3cm}
\begin{remark}
\label{rmk:5}Theorem \ref{theo:3} shows that the correction may be applied for any other scheme satisfying DCC. 
\end{remark}

\section{Convergence and order of accuracy}
\label{sec:Convergence}
 The aim of this section is to prove the following theorem:
\begin{theorem}
\label{theo:3b} Let $u\in C^{2p+3}\left([0,T],X  \right)$ be the exact solution of the problem (\ref{a24}). Suppose that one of the four conditions (i)-(iv) of Theorem \ref{theo:3} holds, with condition (ii) or (iii) holding for all $j=0, 1, \cdots, p+1$. Then each sequence $\left\lbrace u^{2j,n}\right\rbrace _{n=0}^N$, $j=1, 2, \cdots, p+1$, solution of the scheme (\ref{a26}) or  (\ref{a27})-(\ref{a28}),  approximates $u$ with order $2j$ of accuracy. Furthermore, we have the estimate
\begin{equation}
\label{a55}
 \Vert (D_+D_-)^{m}D(u^{2j,n+1/2}-u(t_{n+1/2}))\Vert+\Vert (D_+D_-)^m(u^{2j,n+1}-u(t_{n+1}))\Vert \leq C k^{2j} 
\end{equation} for $m=0,1,...,p-j$ and $n=m+j-1,m+j,..., N-j-m$, where $C$ is a constant depending only on $p$, $T$, and the derivatives of $u$ and $F$ up to order $2m+2j+1$ and $2m+2j-1$, respectively.
\end{theorem}

To prove this theorem we need Theorem \ref{theo:3} and the the following lemma:

\begin{lemma}
\label{lem:2}
Let $\left\lbrace u^{2,n}\right\rbrace _{n=0}^N$ be the solution of the scheme (\ref{a26}). Suppose that one of the conditions (i), (iii) or (iv) of Theorem \ref{theo:3} holds, or $\left\lbrace u^{2,n}\right\rbrace _{n=0}^N$ is bounded in the sense of the condition (ii) of this theorem. Then $\left\lbrace u^{2,n}\right\rbrace _{n=0}^N$ approximates $u$ with order 2 of accuracy, and we have the inequality
\begin{equation}
\label{a43}
 \Vert (D_+D_-)^{m}D(u^{2,n+1/2}-u(t_{n+1/2}))\Vert+\Vert (D_+D_-)^m(u^{2,n+1}-u(t_{n+1}))\Vert \leq C k^{2}, 
\end{equation} for $m=0,1,...,p$ and $n=m,m+1,..., N-m-1$, where $C$ is a constant depending only on $p$, $T$, and the derivatives of $u$ and $F$ up to order $2m+3$ and $2m+1$, respectively.
\end{lemma}

\begin{proof}[Proof of Lemma \ref{lem:2}] For the sake of simplification we suppose that $F=F(x)$. The general case can be handled by transforming (\ref{a24}) to an autonomous system. From the hypotheses of the Lemma, Theorem \ref{theo:3} implies that $\left\lbrace u^{2,n}\right\rbrace _{n=0}^N$ approximates $u$ with order two of accuracy:
\begin{equation}
\label{a42}\Vert u(t_n)-u^{2,n} \Vert \leq Ck^2, \mbox{ for each } n=0,1,2,\cdots, N,
\end{equation}where $C$ is a constant depending only on $T$, $F$ and the derivatives of $u$ up to order 3.
To establish (\ref{a43}) we proceed by induction on the integer $m=0,1,\cdots,p$.
\newline

\item[1.] Inequality (\ref{a43}) for $m=0$.

As in Theorem \ref{theo:3}, we combine (\ref{a24}) and (\ref{a26}) and deduce the identity
\begin{equation}
\label{a44}
D\Theta^{2,n+1/2}=\left[ F\left( \widehat{u}^{2,n+1}\right)-F\left( \widehat{u}(t_{n+1})\right)\right] +\sigma^{2,n+1/2},
\end{equation}where 
 $$\Theta^{2,n}=u^{2,n}-u(t_n),$$
 and 
\begin{align*}
\sigma^{2,n+1/2}=\left[ u'(t_{n+1/2})-Du(t_{n+1/2})\right]  -\left[ F\left( u(t_{n+1/2}\right) -F\left(\widehat{u}(t_{n+1})\right)\right].
\end{align*}
From Taylor's formula with integral remainder and the estimate (\ref{b3}), there exists a function $g$ such that 
$$\sigma^{2,n+1/2} =k^2g(t_{n+1}),$$
with 
\begin{equation}
\label{a48}\Vert  D_+^{m_1}D_-^{m_2} g(t_{n+1})\Vert \leq C, ~\mbox{ for }~ m_2-1\leq n\leq N-m_1-1,
\end{equation}
for each nonnegative integers $m_1$ and $m_2$ such that $m_1+m_2\leq 2p$, where $C$ is a constant depending only on $T$, $F$, and the derivatives of $u$ up to order $m_1+m_2+3$. We can write
$$F\left( \widehat{u}^{2,n+1}\right)-F\left(\widehat{u}(t_{n+1})\right)=\int_0^1dF\left(K^{n+1}_1\right)(\widehat{\Theta}^{2,n+1})d\tau_1,$$
where
$$K^{n+1}_1=\widehat{u}(t_{n+1})+\tau_1 \widehat{\Theta}^{2,n+1}.$$
The last identities substituted into (\ref{a44}) yield 
\begin{equation}
\label{a45}
D\Theta^{2,n+1/2}=\int_0^1dF\left( K^{n+1}_1\right)(\widehat{\Theta}^{2,n+1})d\tau_1+k^2g(t_{n+1}).
\end{equation}
Proceeding as in Theorem \ref{theo:3}, we deduce from (\ref{a42}) and the regularity of $u$ that 
 \begin{equation*}
  \left \Vert \int_0^1dF\left(K^{n+1}_1\right)(\widehat{\Theta}^{2,n+1})d\tau_1\right \Vert \leq C\Vert\widehat{\Theta}^{2,n+1}\Vert.
 \end{equation*} Therefore, taking the norm on both sides of (\ref{a45}), we deduce by the triangle inequality and the inequalities (\ref{a42}) and (\ref{a48}), for $m_1=m_2=0$, that
\begin{equation}
\label{a46}
\Vert D\Theta^{2,n+1/2}\Vert \leq  C\Vert\widehat{\Theta}^{2,n+1}\Vert+k^2\Vert g(t_{n+1})\Vert \leq Ck^2,
\end{equation}
 where $C$ is a constant depending only on $T$ and the derivatives of $u$ and $F$ up to order 3 and 1, respectively. The last inequality combined with (\ref{a42}) implies that (\ref{a43}) holds for $m=0$. 
\newline

\item[2.]  Here we are going to prove that inequality (\ref{a43}) remains true for $m+1$, assuming that it holds for an arbitrary integer $m$ such that $0\leq m\leq p-1$.

We apply $\left( D_+D_-\right)^{m}D_+$ to (\ref{a45}) and obtain 
\begin{equation}
\label{a47} 
\left( D_+D_-\right)^{m+1}\Theta^{2,n+1}=\left( D_+D_-\right)^{m}D_+ h(t_{n+1})  +k^2\left( D_+D_-\right)^{m}D_+g(t_{n+1}),
\end{equation} where we set
$$h(t_{n+1})=\int_0^1dF\left(K^{n+1}_1\right)(\widehat{\Theta}^{2,n+1})d\tau_1.$$
The main difficulty is to bound $\left( D_+D_-\right)^{m}D_+ h(t_{n+1})=D_+^{2m+1}h(t_{n+1-m})$. We have 
$$D_+h(t_{n})=\int_0^1dF\left(K^{n+1}_1\right)(D_+\widehat{\Theta}^{2,n})d\tau_1 +\int_0^1\int_0^1d^2F\left(K^{n}_2\right)\left( D_+K_1^{n},\widehat{\Theta}^{2,n}\right) d\tau_1d\tau_2,$$

\begin{equation*}
\begin{aligned}
&D^2_+h(t_{n})=\int_0^1dF(K^{n+2}_1)(D^2_+\widehat{\Theta}^{2,n})d\tau_1+\int_0^1\int_0^1d^2F(K^{n+1}_2)(  D_+K_1^{n+1},D_+\widehat{\Theta}^{2,n})d\tau^2\\&+\int_0^1\int_0^1d^2F(K^{n+1}_2)( D_+^2K_1^{n},\widehat{\Theta}^{2,n+1})d\tau^2+\int_0^1\int_0^1d^2F(K^{n+1}_2)( D_+K_1^{n},D_+\widehat{\Theta}^{2,n})d\tau^2\\&+\int_0^1\int_0^1\int_0^1d^3F\left(K^{n}_3\right)\left(D_+K_2^n, D_+K_1^{n},\widehat{\Theta}^{2,n}\right) d\tau^3,
\end{aligned}
\end{equation*}
where $d\tau^i =d\tau_1\cdots d\tau_i$, and
\begin{equation}
\label{a47b}
K_{i+1}^{n}=K^n_{i}+\tau_{i+1}(K_{i}^{n+1}-K_{i}^{n})=K_1^n+\sum_{l=1}^i \sum_{2\leq i_1<\cdots <i_l \leq i+1}\tau_{i_1}\cdots\tau_{i_l}k^lD_+^lK_1^n. 
\end{equation}
It follows the general formula
\begin{equation}
\label{a49}
\begin{split}
D_+^qh(t_n)=\sum_{i=1}^{q+1}\sum_{|\alpha_i|=q} L^{n,q}_{i,\alpha_i} ,\mbox{ for } q=1,2,...,2p+1, \mbox{ and }n\leq N-q,
\end{split}
\end{equation}
where $\alpha_i=(\alpha_i^1,\cdots,\alpha_i^{i-1},\alpha_i^i) \in \left\lbrace 1,2,\cdots, q\right\rbrace ^{i-1}\times \left\lbrace 0,1,\cdots, q-i+1\right\rbrace $, and $L^{n,q}_{i,\alpha_i}$ is a linear combination, with properly chosen coefficients, of the quantities
\begin{equation*}
L^{n,q}_{i,\alpha_i,\beta_i} =\int_{[0,1]^i}d^iF(K_i^{n+q+1-i})\left( D_+^{\alpha_i^{i-1}}K_{i-1}^{n+\beta_{i}^{i-1}},\cdots,D_+^{\alpha^1_{i}}K_{1}^{n+\beta_i^1},D_+^{\alpha_i^i}\widehat{\Theta}^{2,n+\beta_i^i}\right)d\tau^i,
\end{equation*}
where $\beta_i=(\beta_i^1,\cdots,\beta_i^{i-1},\beta_i^i) \in \left\lbrace 1,2,\cdots, q\right\rbrace ^{i-1}\times \left\lbrace 0,1,\cdots, q-i+1\right\rbrace $ with $\beta_i^l+\alpha_i^{l}\leq q-l+1$, for $l=1, \cdots, i$. From (\ref{a47b}) and (\ref{a42}) we have 
$$K_i^{n+1}=u(t_{n+1/2})+O(k), \mbox{ for } i=1,2,\cdots,2p+2, $$ and we deduce that there exists $k_3>0$ such that $0<k\leq k_3$ implies
\begin{equation}
\label{a50}\left \Vert d^i F\left( K_{i}^{n} \right) \right \Vert \leq C_i,\;\mbox{ for } i=1,2,...,2p+2, \mbox{ and }n=0,1,\cdots,N-i+1,
\end{equation}
where $C_i$ is a constant depending only on $k_3$, $T$, and the derivatives of $u$ and $F$ up to order $3$ and $i$, respectively. From the inductions hypothesis (\ref{a43}) and inequality (\ref{b3}) we have
\begin{equation}
\label{a51}\Vert D^r_+K_{i}^{n}\Vert \leq C, \mbox{ for } 1\leq r\leq i\leq 2m+3,1\leq n\leq N-i-r+1,
\end{equation}
and 
\begin{equation}
\label{a52}
\Vert D_+^r \widehat{\Theta}^{2,n}\Vert \leq Ck^2, \mbox{ for } 1\leq r\leq 2m+1,1\leq n\leq N-r,
\end{equation}
where $C$ is a constant depending only on $m$, $T$,  and the derivatives of $u$ and $F$ up to order $r+2$ and $r$, respectively. Each $L^{n,q}_{i,\alpha_i,\beta_i}$ being multilinear continuous, we deduce from (\ref{a50})-(\ref{a52}) and the relation $\beta_i^l+\alpha_i^{l}\leq q-l+1$, for $l=1,\cdots,i$, that
$$\Vert L^{n,q}_{i,\alpha_i,\beta_i}\Vert \leq Ck^2 ,\mbox{ for } 1\leq i\leq q+1\leq 2m+2, \, n\leq N-q.$$
It follows by the triangle inequality that (\ref{a49}) for $q=2m+1$ yields 
$$\Vert \left( D_+D_-\right)^{m}D_+ h(t_{n+1}) \Vert =\left \Vert D_+^{2m+1}h(t_{n+1-m})\right \Vert \leq Ck^2,$$
for $n=m,m+1,\cdots, N-(m+1)-1$, where $C$ is a constant depending only on $p$, $T$,  and the derivatives of $u$ and $F$ up to order $2m+4$ and $2m+2$, respectively . Passing to the norm in identity (\ref{a47}), we deduce from (\ref{a48}) and the last inequality that
\begin{equation}
\label{a53} \Vert \left( D_+D_-\right)^{m+1}\Theta^{2,n+1} \Vert \leq Ck^2.
\end{equation} 
Otherwise, applying $D_-$ to (\ref{a47}), inequalities (\ref{a50})-(\ref{a52}) and (\ref{a53}) yield
$$\Vert \left( D_+D_-\right)^{m+1} h(t_{n+1}) \Vert =\left \Vert D_+^{2m+2}h(t_{n-m})\right \Vert \leq Ck^2,$$
for $n=m,m+1,\cdots, N-(m+1)-1$, where $C$ is a constant depending only on $p$, $T$,  and the derivatives of $u$ and $F$ up to order $2m+5$ and $2m+3$, respectively. Therefore, passing to the norm in the identity obtained by applying $D_-$ to (\ref{a47}), we deduce from (\ref{a47}) and the last inequality that 
\begin{equation}
\label{a54} \Vert D_-\left( D_+D_-\right)^{m+1}\Theta^{2,n+1} \Vert \leq Ck^2,
\end{equation}for $n=m,m+1,\cdots, N-(m+1)-1$, with the constant $C$ depending only on $p$, $T$,  and the derivatives of $u$ and $F$ up to order $2m+5$ and $2m+3$, respectively. Inequalities (\ref{a53}) and (\ref{a54}) imply that the induction hypothesis is also true for $m+1$, and we deduce that (\ref{a43}) is true for each integer $m=0,1,...,p$.
\end{proof}

\begin{proof}[Proof of Theorem \ref{theo:3b}] We proceed by induction on $j=1,2,...,p+1$. The case $j=1$ is immediate from Lemma \ref{lem:2}. Suppose that $\left\lbrace u^{2j,n}\right\rbrace _n^N$ approximates $u$ with order $2j$ of accuracy and satisfies (\ref{a55}), for an arbitrary $j$ such that $j\leq p$. We are going to prove that $\left\lbrace u^{2j+2,n}\right\rbrace _n^N$ approximates $u$ with order $2j+2$ of accuracy and (\ref{a55}) holds substituting $j$ by $j+1$. 

From the induction hypothesis, $\left\lbrace u^{2j,n}\right\rbrace _{n}$ satisfies DCC. Because $\left\lbrace u^{2j,n}\right\rbrace _{n}$ and $\left\lbrace \overline{u}^{2j,m}\right\rbrace _m$ are computed from the same scheme DC2j, but for different time steps, $\left\lbrace \overline{u}^{2j,m}\right\rbrace _m$  also satisfies DCC. Therefore, as in \ref{a39a}, Theorem \ref{theo:3} applied to the approximation $\left\lbrace u^{2j+2,n}\right\rbrace _{n=0}^j$, built from $\left\lbrace \overline{u}^{2j,m}\right\rbrace _m$, yields 
\begin{equation*}
\Vert \overline{\Theta}^{2j+2,n} \Vert \leq C\frac{1}{2-\mu k}\left( \frac{2+\mu k}{2-\mu k}\right)^{n-1} k^{2j+2}+\left( \frac{2+\mu k}{2-\mu k}\right)^{n}\Vert \overline{\Theta}^{2j+2,0} \Vert,
\end{equation*}where
$$\overline{\Theta}^{2j+2,n}=\left( u^{2j+2,n}-u(t_{n})\right) -\Gamma^j \left( \overline{u}^{2j,(2j+1)n+j}-u(t_{(2j+1)n+j})\right), \mbox{ for } 1\leq n\leq j.$$
According to the DCC and the condition $u^{2j+2,0}=u(t_0)=u_0$, we have
$$\left \Vert \overline{\Theta}^{2j+2,0} \right \Vert =\left \Vert \Gamma^j \left( \overline{u}^{2j,j}-u(t_{j})\right) \right \Vert \leq Ck^{2j+2}.$$
By the triangle inequality and the DCC, the last two inequalities yield 
\begin{equation}
\label{a54b}
\Vert u^{2j+2,n}-u(t_n) \Vert \leq Ck^{2j+2}, \mbox{ for } n=0,1,\cdots,j.
\end{equation}
From the DCC on $\left\lbrace u^{2j,n}\right\rbrace _n$ and the inequality (\ref{a54b}), Theorem \ref{theo:3} again implies that $\left\lbrace u^{2j+2,n}\right\rbrace _{n=0}^N$ approximates the exact solution $u$ with order $2j+2$ of accuracy. Therefore, it is enough to establish (\ref{a55}) for $j+1$, $j\leq p$. To this end we rewrite identity (\ref{a33}) as follows
\begin{equation}
\label{a56}
\begin{split}
D\Theta^{2j+2,n+1/2}=H(t_{n+1})+\sigma^{2j+2,n+1/2}+ (\Lambda^j -\Gamma^j)D(u^{2j,n+1/2}-u(t_{n+1/2})),
\end{split}
\end{equation}
with
$$H(t_{n+1})=\int_0^1d_uF \left( t_{n+1/2},\widehat{u}(t_{n+1})-\Gamma^j \widehat{u}(t_{n+1})+\tau_1 \widehat{\Theta}^{2j+2,n+1}\right)\left(\widehat{\Theta}^{2j+2,n+1}  \right) d\tau_1,$$
where $\Theta^{2j+2,n}$ and $\sigma^{2j+2,n+1/2}$ are as in Theorem \ref{theo:3}. 
Proceeding as in Lemma \ref{lem:2} and taking the finite difference formulae (\ref{bb3}) and (\ref{bb4}) into account, we can write
\begin{equation*}\sigma^{2j+2,n+1/2} =k^{2j+2}\varepsilon_1 (t_{n+1}),
\end{equation*}
where 
$$\Vert  D_+^{m_1}D_-^{m_2} \varepsilon_1(t_{n+1})\Vert \leq C, \mbox{ for } m_1+m_2\leq 2p-2j \mbox{ and } m_2-1\leq n\leq N-m_1-1,$$
$C$ is a constant depending only on $p$, $T$, and the derivatives of $u$ and $F$. According to the inequality (\ref{a55}) from the induction hypothesis, we may write
$$(\Lambda^j -\Gamma^j)D(u^{2j,n+1/2}-u(t_{n+1/2}))=k^{2j+2}\varepsilon_2(t_{n+1}),$$
where
\begin{equation*}
 \Vert  D_+^{m_1}D_-^{m_2}\varepsilon_2(t_{n+1})\Vert \leq C , \mbox{ for } m_1+m_2\leq 2p-2j+2 \mbox{ and } m_2-1\leq n\leq N-m_1-1.
\end{equation*} Therefore, writing (\ref{a56}) as follows
$$D_-\Theta^{2j+2,n+1}=H(t_{n+1})+k^{2j+2}G(t_{n+1}),$$with 
$$G(t_{n+1})=\varepsilon_1 (t_{n+1})+\varepsilon_2 (t_{n+1}),$$ the induction hypothesis and the reasoning from Lemma \ref{lem:2}, substituting the functions $h$ and $g$, respectively, by $H$ and $G$, $\widehat{\Theta}^{2,n+1}$ by $\widehat{\Theta}^{2j+2,n+1}$, and $k^2$ by $k^{2j+2}$, yields
\begin{equation*}
 \Vert (D_+D_-)^{m}D\widehat{\Theta}^{2j+2,n+1/2}\Vert+\Vert (D_+D_-)^m\widehat{\Theta}^{2j+2,n+1}\Vert \leq C k^{2j+2},
\end{equation*} for $m=0,1,...,p-j$ and $n=m+j-1,m+j,..., N-j-m$, where $C$ is a constant depending only on $p$, $T$, and the derivatives of $u$ and $F$ up to order $2(m+j+1)+1$ and $2(m+j)+1$, respectively. Inequality (\ref{a55}) holds for $\left\lbrace u^{2j+2,n}\right\rbrace_n $ by the triangle inequality from the last inequality.
\end{proof}

We end this section by the following corollary that gives an important convergence property of the DC method. This property is useful for a time-stepping method to solve stiff and large dimensional differential equations arising from the space discretization of time-dependent PDEs.

\begin{corollary} 
\label{cor:saint01}
Suppose that the function $F$ is from $\mathbb{R}^s\rightarrow \mathbb{R}^s$, for a positive integer $s$, and satisfies the one-sided Lipschitz condition (\ref{a31}). Then, each approximate solution $\left\lbrace u^{2j,n}\right\rbrace _{n=0}^N$ from $DC2j$ satisfies the inequality
\begin{equation}
\label{saint2}
\vert u^{2j,n}-u(t_n)\vert\leq Ck^{2j}, \mbox{ for each }  k\in (0,k_0),
\end{equation}
where $C$ is a constant independent from any global Lipschitz constant on $F$, and either $k_0=2/\beta$ for $\beta >0$ or $k_0=+\infty$ for $\beta\leq 0$.
\end{corollary}

\begin{proof} 
From the regularity assumption on $F$ and $u$ and the one sided-Lipschitz condition, we deduce from Theorem \ref{theo:3b} that each $\left\lbrace u^{2j,n}\right\rbrace _{n=0}^N$, $j=1,2,\cdots$, satisfies DCC. Therefore, inequality (\ref{saint2}) is immediate from the part (4) of Theorem \ref{theo:3}. The constant $C$ depends only on the derivatives of $u$ up to order $2j+1$ and, according to (\ref{saint1a})-(\ref{saint1b}) and the mean value theorem, on the bound of the Jacobian $F_y$ on the compact set $[0,T]\times \left\lbrace y\in \mathbb{R}^s : |y|\leq R_j \right\rbrace $.
\end{proof}

\begin{remark} 
\label{rem:saint02}
The convergence property satisfied by the schemes $DC2j$ in Corollary \ref{cor:saint01} is in fact $B$-convergence (see, e.g., \cite{frank1981concept,kraaijevanger1985b}) since the constant $C$ of the global error in (\ref{saint2}) is independent from any global Lipschitz constant of the function $F$. Nevertheless, since in the definition of $B$-\-convergence the constant $C$ depends on high order derivatives of the exact solution $u$, the identity
$$u''(t)=F_t(t,u(t))+F_u(t,u(t))\cdot u'(t)$$
can make any requirement on the independence of the constant $C$ with respect to $F_u$ somewhat artificial. The numerical test on Bernoulli ODE in Section \ref{sec:Numer} gives an application of  Corollary \ref{cor:saint01}.
\end{remark}


\begin{remark}
\label{rem:saint03} In practice, from part 4 of the proof of Theorem \ref{theo:3}, the global error for an approximate solution of the IVP (\ref{a24}) under the one-sided Lipschitz condition (\ref{a31}) by a DC2j+2 method, $j=0,1,2,\cdots$, takes the form
\begin{equation}
\label{a32d}\|u^{2j+2,n}-u(t_n)\|\leq c_{2j+1}C\left( \frac{2+\beta k}{2-\beta k}\right)^{n}k^{2j+2}, ~n\geq j+1,
\end{equation}
for $-2\leq \beta k<2$. The constant $C$ depends on the derivative of the exact solution $u$ of order $2j+3$ and can be very large in magnitude. However, if $\beta<0$ and $k$ is not too small, the factor  $\left( \frac{2+\beta k}{2-\beta k}\right)^{n}$ is sufficiently small such that $C\left( \frac{2+\beta k}{2-\beta k}\right)^{n}<<1$, leading to very accurate approximate solutions for large time steps $k$. Nevertheless, independently of the sign of $\beta$, when $k$ is sufficiently small in the asymptotic region $k\mu<2$, where $\mu$ is the global Lipschitz constant of $F$,  $\left( \frac{2+\beta k}{2-\beta k}\right)^{n}$ becomes closed to 1, for example when $n=j+1$, so that only $c_{2j+3}k^{2j+2}$ must dominate the constant $C$. Consequently, a non B-convergent method can be competitive with respect to a B-convergent one for sufficiently small time steps. This situation will be illustrated by the Bernoulli ODE in Section \ref{sec:Numer}.
\end{remark}

\section{Absolute stability}
\label{sec:Astability}
In this section we prove the absolute stability of the DC schemes. The notion of absolute stability is introduced by Dahlquist \cite{MR0170477} to characterize methods able to solve stiff ODEs. Considering the following IVP,
\begin{equation}
\label{a58c}
\left\lbrace 
\begin{array}{cccc}
u'&=&\lambda u\\
u(0)&=&1,
\end{array}
\right.
\end{equation} where $\lambda$ is a complex number, we have the following definition (see \cite{quarteroni2010,MR0170477}):

\begin{definition}
\label{defn:2}A numerical method is said to be absolutely stable if the corresponding solution for the problem (\ref{a58c}) for fixed $k>0$ and some $Re(\lambda)<0$ is such that 
\begin{equation}
\label{a59}\lim_{n\rightarrow +\infty}|u^n|=0.
\end{equation}
The region of absolute stability of a numerical method is defined as the subset of the complex plane
\begin{equation}
\label{a60}\mathcal{A}=\left\lbrace z=\lambda k \in \mathbb{C}: (\ref{a59}) \mbox{ is satisfied }  \right\rbrace .
\end{equation}
If $\mathcal{A}\cap \mathbb{C}_-=\mathbb{C}_-$, $\mathbb{C}_-=\left\lbrace \lambda \in \mathbb{C} : Re(\lambda)<0 \right\rbrace $, the numerical method is said to be A-stable.
\end{definition}

Before establishing absolute stability results for the deferred correction schemes (\ref{a26}) and (\ref{a27})-(\ref{a28}), we recall the following result.

\begin{lemma}[{ see \cite[formula (6)]{tuenter2006frobenius}  } ]
\label{lem:3}
Let $P_m$ be a polynomial of degree $m$ in one variable. Then the sum $\sum_{i=0}^nP_m(i)$ is a polynomial of degree $m+1$ in the variable $n$.
\end{lemma}

\begin{lemma}
\label{lem:4}Suppose that $F(t,u)=\lambda u$ and $u_0=1$ in the initial value problem (\ref{a24}), where $\lambda$ is a complex number with negative real part ($\lambda \in \mathbb{C}_-$). Then the corresponding approximate solutions from the schemes (\ref{a26}) and (\ref{a27})-(\ref{a28}) can be written as follows
\begin{equation}
\label{a61} u^{2j+2,n}=\left( \frac{2+\lambda k}{2-\lambda k}\right)^{n-j}P_j\left( n\right) ,\mbox{ for } j=0,1,2,..., \mbox{ and } n\geq j,
\end{equation} where $P_j(n)$ is a polynomial of degree $j$  in the variable $n$. 
\end{lemma}

\begin{proof}We suppose that $\lambda k \neq -2$, otherwise we trivially have $u^{2j,n+1}=0$, for $n\geq j$. Since $F(t,u)=\lambda u$, we can rewrite (\ref{a27}) as follows
$$u^{2j+2,n+1}=\frac{2+\lambda k}{2-\lambda k}u^{2j+2,n}+\frac{2}{2-\lambda k}\left( kD_-\Lambda^j u^{2j,n+1}-\lambda k\Gamma^j \widehat{u}^{2j,n+1} \right)$$
where, according to formulae (\ref{bb1}) and (\ref{bb2}), we have
\begin{align*}
kD_-\Lambda^j u^{2j,n}&=\sum_{i=1}^jc_{2i+1}k^{2i+1}D_-(D_+D_-)^iu^{2j,n}\\&=\sum_{i=1}^j\sum_{m=0}^{2i+1}c_{2i+1}(-1)^m\binom{2i+1}{m}u^{2j,n+i-m},
\end{align*}and
$$\Gamma^j \widehat{u}^{2j,n}=\sum_{i=1}^jc_{2i}k^{2i}(D_+D_-)^i\widehat{u}^{2j,n}=\sum_{i=1}^j\sum_{m=0}^{2i}c_{2i}(-1)^m\binom{2i}{m}\widehat{u}^{2j,n+i-m}.$$
Combining the last three identities, we deduce that
\begin{equation}
\label{a62}u^{2j+2,n+1}=\frac{2+\lambda k}{2-\lambda k}u^{2j+2,n}+\frac{2}{2-\lambda k}\sum_{i=0}^{2j+1}\alpha_{j,i}(\lambda k)u^{2j,n+1+j-i},\mbox{ for } n\geq j\geq 1,
\end{equation} where $\alpha_{j,i}$ is affine in $\lambda k$. Under the hypothesis of the lemma, (\ref{a26}) matches the trapezoidal rule, and we have  
\begin{equation*}
u^{2,n}=\left( \frac{2+\lambda k}{2-\lambda k}\right)^{n},
\end{equation*}that is (\ref{a61}) is true for $j=0$. Suppose that (\ref{a61}) holds for an arbitrary integer $j\geq 0$. From (\ref{a62}) we have 
$$u^{2j+4,n}=\frac{2+\lambda k}{2-\lambda k}u^{2j+4,n-1}+\frac{2}{2-\lambda k}\sum_{i=0}^{2j+3}\alpha_{j+1,i}(\lambda k)u^{2j+2,n+1+j-i},$$
with $n\geq j+2$, and, substituting each $u^{2j+2,n+1+j-i}$ by the formula given by the induction hypothesis (\ref{a61}), we deduce that
$$u^{2j+4,n}=\frac{2+\lambda k}{2-\lambda k}u^{2j+4,n-1}+\left( \frac{2+\lambda k}{2-\lambda k}\right)^{n-j-1}Q_{j}(n),$$
 where
$$Q_{j}(n)=\frac{2}{2-\lambda k}\sum_{i=0}^{2j+2}\alpha_{j+1,i}(\lambda k)\left( \frac{2+\lambda k}{2-\lambda k}\right)^{j+2-i}P_{j}(n+1+j-i).$$
It follows that
$$u^{2j+4,n}=\left( \frac{2+\lambda k}{2-\lambda k}\right)^{n-j-1}\left(u^{2j+4,j+1} +\sum_{i=j+2}^{n}Q_{j}(i)\right). $$
It is clear that $Q_{j}(n)$ is a polynomial of degree $j$ in the variable $n$ as $P_j(n)$. Therefore, according to the Lemma \ref{lem:3}, $\sum_{i=j+2}^{n}Q_{j}(i)$ is a polynomial of degree $(j+1)$ in the variable $n$. Whence, 
\begin{equation*}u^{2j+4,n}=\left( \frac{2+\lambda k}{2-\lambda k}\right)^{n-j-1}P_{j+1}(n),\; n\geq j+1,
\end{equation*}
where
$$P_{j+1}(n)=u^{2j+4,j+1} +\sum_{i=j+2}^{n}Q_{j}(i)$$ is a polynomial of degree $j+1$ in the variable $n$. We then deduce by induction that the lemma is true for arbitrary non-negative integer $j$.
\end{proof}
\begin{theorem}
\label{theo:6}Each of the deferred correction schemes (\ref{a26}) and (\ref{a27})-(\ref{a28}) is A-stable.
\end{theorem}

\begin{proof}
From Lemma \ref{lem:4} we have, for $Re(\lambda k)<0$,
$$\lim_{n\rightarrow +\infty}|u^{2j+2,n}|=\lim_{n\rightarrow +\infty}\left |\left( \frac{2+\lambda k}{2-\lambda k}\right)^{n-j}P_j\left( n\right)\right |=\lim_{n\rightarrow +\infty}|P_j\left( n\right)|e^{(n-j)ln \left |\frac{2+\lambda k}{2-\lambda k}\right |}=0$$
since, under the condition $Re(\lambda k)<0$, we have $\left |\frac{2+\lambda k}{2-\lambda k}\right |<1$.
\end{proof}

\section{Numerical experiments}
\label{sec:Numer}

In this section we evaluate the accuracy and order of convergence of the schemes $DC2, DC4, \cdots , DC10$, implemented using the Scilab programming language. The starting values are computed using the scheme (\ref{a73})-(\ref{a74b}).

We choose six standard problems for the evaluation. The first problem concerns $B$-\-convergence by considering a Bernoulli equation. The second problem is about long term integration with an oscillatory solution of large amplitude. The four other problems are about stiffness. The third and fourth problems (B5 modified and E5, respectively) both involve complex eigenvalues of negative real parts, where the imaginary parts of the eigenvalues for the third problem have larger magnitudes while those from the fourth problem have smaller magnitudes. The fifth problem (Robertson) is nonlinear and stiff with real negative eigenvalues, and it also addresses B-convergence. The sixth problem is the van der Pol oscillator, which is stiff with arbitrary complex eigenvalues.

The first three problems have analytic solutions. For problems (\ref{a72}), (\ref{69}) and (\ref{a71}) that do not have an analytic solution, we consider a small time step such that the approximate solutions with $DC6, \cdots, DC10$ are almost identical (to machine precision for problem (\ref{69})), and we choose one of the approximate solutions as reference solution.

For solutions $u=(u_1,\cdots,u_d)~:~[0,T] \rightarrow \mathbb{R}^d$, $1\leq d\leq 6$, the absolute error on the approximate solutions $\left\lbrace u^{2j,n} \right\rbrace_{0\leq n\leq N} $, $1\leq j\leq 5$, is computed with the norm
$$\|u^{2j}_i-u_i\| =\max _{0\leq n\leq N}|u^{2j,n}_i-u_i(t_n)|,\quad 1\leq i\leq d.$$ 
For very large $N$ we extract solutions at $2\times 10^6$ or $3\times 10^6$ discrete times evenly spread over the interval $[0,T]$.

For a comparison of accuracy, we implement in Scilab the backward differentiation formulae (BDF) of order 2, 4 and 6, and the explicit Runge-Kutta (RK) of order 4. The implemented BDF are run with exact starting values for the first three problems that have analytic solutions, while for problems four and five the starting values are provided by the function \texttt{stiff} (implementing BDF with adaptive steps) of the solver \texttt{ode} from Scilab. For the van der Pol oscillator, the comparison of our DC methods is done only with the solutions from \texttt{stiff} and \texttt{rkf} from the solver \texttt{ode}. For each of the problems, except the first one, we give a table of absolute errors and orders of convergence for pairs of two consecutive time steps, for the approximate solutions with the DC methods. We denote by $k_{max}$ the maximal time step allowed to compute an approximate solution with the solver \texttt{stiff} or \texttt{rkf} (see \cite{enright1975comparing} for a discussion on maximal time steps).

\subsection{Bernoulli differential equation}
\begin{equation}
\label{b77b}
u'(t)=F(t,u)=-0.1u(t)-1000u^{20}(t)~,~~~u(0)=1, ~~t\in [0,10].
\end{equation}
 Table \ref{tab:4a}  gives the absolute error and the order of convergence for each pair of consecutive time steps, in the case of DC, BDF and RK4 methods. The dash for RK4 indicates that the method is unstable for the corresponding time steps.

\begin{table}[!ht]
\caption{Absolute error (order of convergence) for the Bernoulli problem.}
\label{tab:4a}       
\begin{tabular}{llllllll}
\hline\noalign{\smallskip}
 $k$ &  \centering DC2 & \centering DC4 & \centering DC6 &\centering DC8& DC10\\[.35ex]
 1 & 0.18 & 1.7e-2 & 1.8e-4 & 2.3e-4 & 1.3e-4  \\[0.35ex]
 2.03e-3 & 3.71e-2~(0.26)   & 6.16e-4~(0.53)  & 7.14e-5~(0.14) & 1.47e-6~(0.81) & 9.42e-7~(0.79)\\[.35ex]
 1.00e-4 & 1.92e-3~(0.98) & 2.93e-5~(1.01) & 4.31e-6~(0.93) & 3.72e-7~(0.45) &5.78e-8~(0.94) \\[.35ex]
 1.00e-5 & 2.22e-5~(1.94) & 1.30e-7~(2.35) & 3.92e-9~(3.04) & 1.9e-10~(3.27) &1.1e-11~(3.73) \\[.35ex]
 5.00e-6 & 5.55e-6~(2.0) & 1.04e-8~(3.70) & 1.4e-10~(3.70) & 4.4e-12~(5.50) & 4.4e-13~(4.64) \\[.35ex]
 3.33e-6 & 2.46e-6~(1.99) & 2.59e-9~(3.33) & 1.6e-11~(5.31) & 4.5e-13~(5.63) & 2.0e-13~(2.02) \\[.35ex]
 2.25e-6 & 1.39e-6~(1.99) & 8.7e-10~(3.79) & 3.3e-12~(5.54) & 4.2e-13~(0.16) & 4.2e-13~(-2.66) \\[.35ex]
\hline\noalign{\smallskip}
 $k$ &  \centering BDF2 & \centering BDF4 & \centering BDF6 & \centering \texttt{RK4} & \\[.5ex]
 1       & 0.14   & 0.83 & 6.1e-2 & -- &  \\[.35ex]
 2.03e-3 & 4.3e-2~(0.19) & 2.5e-2~(0.19) & 1.9e-3~(0.19) & -- &  \\[.35ex]
1.00e-4     & 6.61e-3~(0.62) & 2.98e-3~(0.71) & 1.79e-3~(0.79) & 1.27e-3~ &  \\[.35ex]
1.00e-5     & 2.59e-4~(1.41) & 1.92e-5~(2.19) & 3.15e-6~(2.76) & 4.91e-8~(4.41) &  \\[.35ex]
5.00e-6     & 7.29e-5~(1.91) & 1.92e-6~(3.58) & 1.35e-7~(5.11) & 2.53e-9~(4.28) &  \\[.35ex]
\hline\noalign{\smallskip}
\end{tabular}
\end{table}

 This problem addresses $B$-\-convergence since the function $F$ is one-sided Lipschitz with $\beta=-0.1$, when positive solutions are considered. Moreover,  the problem is strongly nonlinear with exponentially increasing magnitude of derivatives of the right side function $F$. Such derivatives of large magnitude generally limit the accuracy of high order methods that are not B-convergent. The one-sided Lipschitz constant being negative, in accordance with Corollary \ref{cor:saint01}, DC methods provide very accurate approximate solutions for large time steps, and their accuracy increases with the order of the method. However, the convergence of the DC methods is suboptimal, due to the effect of the strong nonlinearity of the ODE. While $DC4$ and $DC6$ almost achieve their proper order for $k\leq 3\times 10^{-5}$, the order of convergence of $DC8$ and $DC10$ are not observed since these methods quickly achieve machine accuracy. In fact, DC10 achieves order 8.05 of convergence for $k=7.14\times 10^{-6}$ to $k=6.66\times 10^{-6}$. BDF methods are stable for large time steps, but they are less accurate than  their corresponding DC methods. RK4 is completely unstable for $k\geq 2.03\times 10^{-3}$. For sufficiently small time steps in the asymptotic region, RK4 is more accurate than DC4 and any of the BDF methods, as stated in Remark \ref{rem:saint03}, while DC6-10 achieve better accuracy.

\subsection{Oscillatory problem \cite{hull1972comparing}}
\begin{equation}
\label{66}
u'=\lambda u\cos (t)~,~~~u(0)=1, ~~T=10^6, \lambda=10.
\end{equation}
The exact solution is  $u(t)=e^{\lambda \sin(t)}$. The original problem is set with $\lambda =1$ in \cite{hull1972comparing}. The author in \cite{karouma2015class} solved this problem with Runge-Kutta methods of orders 4 and 8, for $\lambda =2$ and $T=2580\pi$, to ``illustrate the need of higher order methods when a long-term integration problem is considered''. Table  \ref{tab:4}  gives the absolute error and the order of convergence for each pair of consecutive time steps. The BDF methods are run only for the smallest time step. The solvers \texttt{rkf} and \texttt{stiff} use adaptive time stepping with a maximal time step $k_{max}=0.1$ and tolerances $rtol=100\times atol=10^{-10}$.

\begin{table}[!ht]
\caption{Absolute error (order of convergence) for the oscillatory problem.}
\label{tab:4}       
\begin{tabular}{llllllll}
\hline\noalign{\smallskip}
 $k$ &  \centering DC2 & \centering DC4 & \centering DC6 &\centering DC8 & DC10\\[.5ex]
 5.00e-2 & 3418 & 456.26 & 42.665 & 3.2350 & 0.2132  \\[0.5ex]
 2.50e-2 & 790.2~(2.1) & 25.351~(4.2)  & 0.5959~(6.2)  & 1.17e-2~(8.1) &1.9e-4~(10.1)\\[.5ex]
 1.25e-2 & 193.8~(2.0) & 1.5493~(4.0)  & 9.17e-3~(6.0) & 5.28e-5~(7.8) & 2.79e-6~(6.1)\\[.5ex]
 6.25e-3 & 48.23~(2.0) & 9.67e-2~(4.0) & 1.4e-4~(5.99) & 2.78e-6~(0.0) &2.78e-6~(0.0) \\[.5ex]
 1.56e-3 & 3.010~(2.0) & 3.8e-4~(3.99) & 4.72e-6~(2.5) & 4.67e-6~(-0.3) & 4.7e-6~(-0.3) \\[.5ex]
\hline\noalign{\smallskip}
 $k$ &  \centering BDF2 & \centering BDF4 & \centering BDF6 & \centering \texttt{rkf} & \texttt{stiff}\\[.5ex]
 1.56e-3 & 22026.46 & 14836.76 & 5578.40 & 22026.46 & 2636.00 \\[.5ex]
\hline\noalign{\smallskip}
\end{tabular}
\end{table}

The magnitude of the exact solution $u(t)=e^{10\sin(t)}$ of the modified oscillatory problem is large, resulting in a relatively large absolute error obtained by the DC schemes (absolute errors of about $10^{-7}$ is possible for a good choice of stepsize). Moreover, the long term integration influences the accuracy of these schemes since they achieve absolute errors of about $10^{-9}$ when the solution interval is reduced to $[0,1000]$. Nevertheless, each DC scheme converges with its proper order. The DC methods are considerably more accurate than standard methods (both with fixed and variable stepsizes) which are inaccurate for this problem. For instance, for BDF2 and \texttt{rkf}, the solutions remain bounded with bounds close to the maximal amplitude of the exact solution but the phase of the oscillation is completely wrong.

\subsection{Problem B5 modified \cite{enright1975comparing}, stiff with complex eigenvalues of negative real parts and larger (in magnitude) imaginary parts} \label{sec:B5}
\begin{equation}
\label{a70}
y'=\begin{bmatrix}-10 & ~~\alpha & ~~0 & ~~0 & 0 & 0\\
               -\alpha & -10     & ~~0 & ~~0 & 0 & 0\\
                   ~~ 0 & ~~0      & -4&~~ 0 & 0 & 0\\
                   ~~ 0 & ~~0      & ~~0 & -1& 0 & 0\\
                   ~~ 0 & ~~0      & ~~0 &~ ~0 & -0.5 & 0\\
                    ~~0 & ~~0      & ~~0& ~~~0~& 0 &-0.1\end{bmatrix}y,~ y(0)=\begin{bmatrix}
                    1\\1\\1\\1\\1\\1
                    \end{bmatrix},~ \alpha=5000,~ T=20.
\end{equation}
This problem, originally set with $\alpha =100$, is an illustration of ODEs resulting from a semi-discretization by finite element methods of parabolic PDEs \cite{stewart1990avoiding}. We choose $\alpha =5000$ to make the problem a little more difficult. Table \ref{tab:5}  gives the absolute errors for the first component of the approximate solutions which is similar for the second component. The absolute errors for the others components quickly achieve machine precision. The solvers  \texttt{stiff} and \texttt{rkf} are run with $k_{max}=2\times 10^{-5}$ and $atol=10\times rtol=10^{-15}$.

The imaginary parts of the Jacobian eigenvalues of the modified B5 problem are large. Even though the real parts of the eigenvalues are negative, we observe that smaller time steps are required by DC schemes to obtain accurate approximations. DC schemes achieve their proper order of convergence, but BDF methods perform better for this problem than DC schemes.

\begin{table}[!ht]
\caption{Absolute error (order of convergence) for the first component of the solution for $B5$ modified}
\label{tab:5}       
\begin{tabular}{llllllll}
\hline\noalign{\smallskip}
 $k$ &  \centering DC2 & \centering DC4 & \centering DC6 &\centering DC8 & DC10\\
 2.000e-5 & 0.2152 & 6.51e-2 & 2.22e-2 & 8.00e-3 & 2.98e-3  \\
5.000e-6 & 1.35e-2~(2) & 2.59e-4~(4) & 5.59e-6~(6) & 1.27e-7~(8)   & 2.97e-9~(10)\\
2.500e-6 & 3.38e-3~(2) & 1.62e-5~(4) & 8.74e-8~(6) & 4.9e-10~(8)  & 2.9e-12~(10) \\
1.250e-6 & 8.47e-4~(2) & 1.01e-6~(4) & 1.36e-9~(6) & 1.9e-12~(8)   & 7.4e-14~(5.3)\\
3.125e-7 &5.29e-5~(2) & 4.00e-9~(4) & 3.6e-13~(6)  & 7e-14 (2.4)& 6.3e-14\\[.5ex]
6.250e-8  & 2.11e-6~(2) & 6.3e-12~(4) & 6.02e-13     & 2.33e-13~      & 1.19e-13~ \\
\noalign{\smallskip}\hline
 $k$ &  \centering BDF2 & \centering BDF4 & \centering BDF6 & \centering \texttt{rkf} & \texttt{stiff}\\[.5ex]
 1.25e-6 & 3.38e-3 & 7.94e-8 & 2.3e-12 & 2.36e-6 & 6.6e-10 \\[.5ex]
\hline\noalign{\smallskip}
\end{tabular}
\end{table}

\subsection{Problem E5\cite{enright1975comparing}, stiff with complex eigenvalues of negative real parts and smaller (in magnitude) imaginary parts} \label{sec:E5}
\begin{equation}
\label{a72}
 \begin{aligned}
y'_1&=-7.89\times 10^{-10}y_1-1.1\times 10^7y_1y_2&\\
y'_2&=7.89\times 10^{-10}y_1-1.13\times 10^9y_2y_3&\\
y'_3&=7.89\times 10^{-10}y_1-1.1\times 10^7 y_1y_2+1.13\times 10^3y_4-1.13\times 10^9y_2y_3&\\
y'_4&=1.1\times 10^7y_1y_2+1.13\times 10^3y_4&\\
\displaystyle &y(0)=(1.76\times 10^{-3},0;0;0)^t, T=1000.
\end{aligned}
\end{equation}
A reference solution is computed with $DC10$ for $k= 10^{-3}$. The solution of this problem has small magnitude in $[1.618\times 10^{-3},1.76\times 10^{-3}]\times [0,1.46\times 10^{-10}]\times [0,8.27\times 10^{-12}]\times [0,1.38\times 10^{-10}]$ and the eigenvalues of the Jacobian matrix $dF(y)$ along the solution curve belong to the region  
$[-20490,3.68\times 10^{-12}]\times [-9.17\times 10^{-5},9.17\times 10^{-5}]$ of the complex plane. Table \ref{tab:9} gives the absolute errors and order of convergence for the four components of the approximate solutions. For BDF, RK4 and \texttt{stiff}, the absolute errors are provided only for the first component. The  absolute error on the other components is smaller by 2 (RK4) to 5 (\texttt{stiff}) orders of magnitude, as we should expect from the magnitude of the solution components. The implemented BDF methods are run with starting values deduced from the solver \texttt{stiff}. The implemented RK4 is unstable for time steps $k\geq 2\times 10^{-4}$, and the absolute error is reported for $k=10^{-4}$ in table~\ref{tab:9}. The solver \texttt{stiff} is run with $k_{max}=10^{-3}$ and $rtol=10^8\times atol=10^{-15}$.

\begin{table}[!ht]
\caption{ Absolute error (order of convergence) for the problem E5}
\label{tab:9}
\begin{tabular}{lllllllll}
\hline\noalign{\smallskip}~ $k$ &  \centering DC2 & \centering DC4 & \centering DC6 &\centering DC8 & DC10 \\ 
\hline\noalign{\smallskip}
\multirow{4}{4em}{100}
&2.79e-07 & 5.34e-08 & 8.31e-09& 4.26e-09 & 1.04e-09\\
&8.30e-12 & 9.68e-13 & 6.86e-14& 6.14e-14 & 1.66e-14\\
&4.47e-13 & 5.31e-14 & 3.28e-15& 3.40e-15 & 8.42e-16\\
&7.85e-12 & 9.14e-13 & 6.54e-14& 5.81e-14 & 1.57e-14\\
\hline\noalign{\smallskip}
\multirow{4}{4em}{50}
   &7.52e-08(1.89) & 1.02e-08(2.38) & 1.56e-09(2.41) &  8.53e-11(5.64) & 4.92e-11(4.41)\\
   &1.96e-12(2.08) & 6.46e-14(3.90) & 3.16e-14(1.12) &  2.94e-15(4.38) & 5.07e-16(5.03)\\
   &1.07e-13(2.06) & 3.73e-15(3.83) & 1.61e-15(1.02) &  2.21e-16(3.94) & 9.78e-17(3.11)\\
   &1.86e-12(2.08) & 6.14e-14(3.89) & 3.00e-14(1.12) &  2.85D-15(4.35) &  4.09D-16(5.26)\\
\hline\noalign{\smallskip}
\multirow{4}{4em}{10}
& 3.16e-09(1.99) & 2.37e-11(4.03) & 5.26e-13(5.23) & 1.28e-14(6.72) & 4.51e-16(8.89)\\
& 7.77e-14(1.99) & 2.79e-16(3.68) & 3.02e-18(5.74) & 1.15e-19(7.94) & 7.28e-21(8.09)\\
& 4.31e-15(1.97) & 7.08e-17(1.79) & 5.91e-17(0.24) & 6.27e-17(0.12) & 6.84e-17(0.09)\\
& 7.34e-14(1.99) & 3.20e-16(3.37) & 6.18e-17(1.79) & 6.28e-17(0.57) & 6.84e-17(0.11) \\
\noalign{\smallskip}\hline
 $k$ &  \centering BDF2 & \centering BDF4 & \centering BDF6 & \centering RK4 & \texttt{stiff}\\[.5ex]
 10 & 5.7e-8 & 6.6e-10 & 3.5e-11 & 2.03e-16 & 1.29e-16 \\[.5ex]
\hline\noalign{\smallskip}
\end{tabular}
\end{table}

Imaginary parts of eigenvalues for the problem E5 are smaller, and larger time steps allow DC schemes to produce very accurate approximations, compared to the modified B5 problem. DC schemes perform better for this problem than BDF methods. They achieve their proper order of convergence but on a relatively small range of time steps, for higher order DC methods, since the solution is already very accurate for large time steps.

\subsection{Robertson (1966) \cite{wanner1991solving}, stiff with real negative eigenvalues}
\begin{equation}
\label{69} \begin{aligned}
y'_1&=-0.04y_1+10^{4}y_2y_3\\
y'_2&=0.04y_1-10^{4}y_2y_3-3.10^{7}y_2^2\\
y'_3&=3.10^{7}y_2^2\\
y(&0)=(1,0,0)^t,~~T=10^5.
\end{aligned}
\end{equation}
This is one of the three problems considered as stiffest in \cite{wanner1991solving}. We compute a reference solution with DC10 for the time step $k=1/6000$. The solution belongs to the region $[1.78\times 10^{-2},1.00]\times [0,3.58\times 10^{-5}]\times [0,0.983]$ and the eigenvalues of the Jacobian dF(y) along the solution curve belong to $[-9825.744,0]$. Table \ref{tab:6} gives absolute errors and orders of convergence of DC methods for each component of the solution. For other methods, we give only the maximal errors on the three components of the approximate solutions. The solver \texttt{stiff} is run with $k_{max}=1/600$ and $rtol=100\times atol=10^{-15}$. The solver \texttt{rkf} fails in solving this problem for various tolerances and $k_{max}$, and Scilab reported: ``it is likely that rkf45 is inefficient for solving this problem''. The implemented BDF methods are run with starting values deduced from the solver \texttt{stiff} using the preceding tolerances.
 
\begin{table}[!h]
\caption{Absolute error (order of convergence) for Robertson problem }
\label{tab:6}     
\renewcommand{\arraystretch}{1}\begin{tabular}{llllllll}
\hline\noalign{\smallskip}~ $k$ &  \centering DC2 & \centering DC4 & \centering DC6 &\centering DC8 & DC10\\[1.2ex] 
\hline\noalign{\smallskip}
\multirow{3}{3.5em}{0.5}
   &3.63e-5 & 4.46e-6 & 2.08e-6 & 2.91e-6 & 3.09e-6\\
   &3.63e-5 & 4.46e-6 & 2.08e-6 & 2.91e-6 & 3.09e-6\\
   &7.12e-5 & 4.37e-7 & 1.02e-7 & 4.12e-7 & 4.26e-7\\
\hline\noalign{\smallskip}
\multirow{3}{3.5em}{1/300}
   & 4.7e-9~(1.8) & 1.09e-9~(1.7)  &  4.0e-10~(1.7)  & 3.0e-10~(1.9)  & 2.0e-10~(1.9)\\
   & 7.4e-9~(1.7) & 2.23e-8~(1.1)  &  4.16e-8~(0.8)  & 2.9e-8~(0.9) & 2.5e-8~(0.9)\\
   & 4.7e-9~(1.9) & 2.12e-8~(0.6)  &  4.12e-8~(0.6)  & 2.8e-8~(0.5) & 2.5e-8~(0.6)\\
\hline\noalign{\smallskip}
\multirow{3}{3.5em}{1/600} 
& 1.0e-9~(2.2)   & 1.5e-10~(2.8)  & 1.0e-12~(8.6)    & 9.9e-13~(8.) & 7.5e-13~(8.2)\\
& 5e-13~(14.)    & 3.0e-14~(19.6)    & 2.0e-16~(27.7)   & 2.0e-16~(27.1)  & 3.0e-16~(26.1)\\
& 1.0e-9~(2.2)   & 1.5e-10~(7.1)  & 1.0e-12~(15.3)   & 9.9e-13~(15)& 4.0e-13~(15.8)\\
\hline\noalign{\smallskip}
\hline\noalign{\smallskip}
\multirow{3}{3.5em}{1/6000}
&  9.24e-12 &  7.31e-14 &  1.48e-14 &  4.57e-14  & --\\
&  5.38e-15 &  0.       &  0.       &  0.        & --\\
&  9.25e-12 &  2.07e-13 &  1.36e-13 &  8.27e-14  & --\\
\noalign{\smallskip}\hline
\noalign{\smallskip}\hline
 $k$ &  \centering BDF2 & \centering BDF4 & \centering BDF6 & \centering \texttt{RK4} & \texttt{stiff}\\[.5ex]
 0.5 & 5.3e-4 & 3.6e-5 & 4.1e-6 & -- & 7.76e-13 \\[.5ex]
 1/600 & 2.8e-6 & 1.2e-6 & 6.9e-7 & -- & 7.28e-13 \\[.5ex]
\noalign{\smallskip}\hline
\end{tabular}
\end{table} 

The Robertson problem is stiff and addresses B-convergence since its Jacobian matrix has real negative eigenvalues with some having large magnitude. For this problem, DC schemes produce accurate approximate solutions even for large time steps, and high order DC methods can be avoided (DC6 is enough). The convergence is slow for $k>1/300$, but faster convergence happens for $k$ in the asymptotic region ($k<1/300$). The DC schemes perform better than BDF methods at equal order and time step. A comparison of the errors for $k=1/600$ suggests that the error constants might be 3 to 5 orders of magnitude smaller for DC than BDF methods.

\subsection{ van der Pol oscillator \cite{enright1975comparing,shampine1981evaluation}, stiff with arbitrary complex eigenvalues}
\begin{equation}
\label{a71}
 \begin{array}{ccc}
\displaystyle y'_1=y_2 ~~~~~~~~~~~~~~~~~~~~~~~~~~~~~~~~~~~~~~~~~~~~~\\

\displaystyle y'_2=\mu (1-y_1^2)y_2-y_1~~~~~~~~~~~~~~~~~~~~~~~~~~~\\

\displaystyle y_1(0)=2,~~y_2(0)=0, T=3000,\mu=1000.
\end{array}
\end{equation}
This problem was initially proposed for $T=1$ and $\mu=5$ in \cite{enright1975comparing}. The actual version results from a  suggestion by Shampine \cite{shampine1981evaluation}. We compute a reference solution with DC8 for  $k=1.875\times\times 10^{-6}$. The solution belong to the region $[-2, 2.000073]\times [-1323.04,1231.35]$ of the real plane and the eigenvalues along the solution curve belong to the region $[-3000.29,1123.17]\times [-1158.48,1158.48]$ of the complex plane. Table \ref{tab:8}  gives the absolute errors and orders of convergence.
For \texttt{rkf} and \texttt{stiff}, we use $k_{max}=7.5\times 10^{-5}$ and $rtol=10$, $atol=10^{-16}$.

\begin{table}[!h]
\caption{Absolute error (order of convergence) for the van der Pol's equation }
\label{tab:8}     
\renewcommand{\arraystretch}{1}\begin{tabular}{cccccccc}
\hline\noalign{\smallskip}
$k$ &  \centering DC2 & \centering DC4 & \centering DC6 &\centering DC8 & DC10 \\
\hline\noalign{\smallskip}
\multirow{2}{3.65em}{3.75e-5}
&3.0089 & 2.9999 & 2.9440 & 0.1838 & 3.12e-3\\
&1322.9 & 1327.5 & 1320.6 & 197.79 & 3.26792\\
\hline\noalign{\smallskip}
\multirow{2}{3.65em}{1.50e-5}
&2.9769~(0) & 2.9999~(0) & 0.1080~(3.6) & 1.90e-4~(7.5) & 5.1e-5~(4.5)\\
&1333.3~(0) & 1330.3~(0) & 113.69~(2.7) & 0.18281~(7.6) & 5.1e-2~(4.5)\\
\hline\noalign{\smallskip}
\multirow{2}{3.65em}{7.50e-6}
&2.8706~(0) & 2.6947~(0) & 1.60e-3~(6.0) & 1.74e-6~(6.7) & 1.27e-5~(1.9)\\
&1327.4~(0) & 1286.5~(0) & 1.6349~(6.1) & 1.80e-3~(6.7) & 1.29e-2~(1.9)\\
\hline\noalign{\smallskip}
\multirow{2}{3.65em}{1.875e-6}
&0.74(0.9) & 0.339~(1.5) & 2.50e-7~(6.3) & -- & 2.88e-7~(2.7)\\[.2ex]
&659.~(0.5)& 373.2~(0.9) & 2.91e-4~(6.2) & -- & 2.92e-4~(2.7)\\[.2ex]
\noalign{\smallskip}\hline
-- & \centering \texttt{stiff} & \centering \texttt{rkf} &  &  &  \\
\hline\noalign{\smallskip}
\multirow{2}{3.65em}{}
& 2.16e-6 & 3.54e-2 &  &  & \\
& 3.48e-3 & 64.76 &   &  & \\
\hline\noalign{\smallskip}
\end{tabular}
\end{table}

The van der Pol oscillator is stiff and the solution has a large magnitude. DC6 and DC8 reached their order of convergence. This shows that the DC strategy works well in spite of the fact that DC2 and DC4 would require much smaller time steps to produce reasonably accurate solutions. The order of convergence for DC10 is not observed, though the solutions obtained are accurate.

\subsection{Discussion of the numerical results} 

 In general, a careful assessment of the proof of Theorem \ref{theo:3} points out to the fact that, for a system with complex eigenvalues $\lambda =\lambda_1+i\lambda_2$, we only need a time step $k$ such that $k\,\max\left\lbrace \lambda_1,|\lambda_2|\right\rbrace <2$ for a good accuracy (faster convergence happens when $- \lambda_1>>|\lambda_2|$). These situations are well illustrated by the test cases of Sections \ref{sec:B5} and \ref{sec:E5}, where the required time step for accuracy is much smaller for modified B5 than E5. However, time steps $k$ such that $k\mu\simeq k|\lambda|<2$, $\mu\simeq\displaystyle \max_{0\leq t\leq T}\Vert d_uF\left( t,u(t)\right)\Vert$, is necessary for an asymptotic convergence with proper order. For example, in the case of the Bernoulli equation we have $\lambda \simeq -20000.1<0$ and $\mu =20000.1$. Large time steps provide accurate approximations (as expected from B-convergent methods), but asymptotic convergences are observed only for $k\mu<2$.

For the computational effort of the DC methods, we recall that to compute an approximate solution on discrete points $0=t_0<t_1<\cdots <t_N=T$, $DC2$ solves $N$ nonlinear systems while $DC2j$, $j\geq 2$, solves $j\times N$ systems. In the case of the Bernoulli equation, for example, $DC10$ achieves the maximal error of about $1.1\times 10^{-11}$ by solving approximately $5\times 10^6$ nonlinear systems while the maximal absolute error for $DC2$ is about $8.9\times 10^{-7}$ for $N=5\times 10^6$. We did not report any CPU time since our code is written in Scilab, an interpreted language. All methods that we implemented are consequently interpreted, while \texttt{rkf} and \texttt{stiff} provided with Scilab are compiled. Nevertheless, the main burden in implicit time-stepping solvers is the resolution of nonlinear systems, and we have shown that higher order DC methods give the most accurate approximations by solving fewer systems of equations. This gives a clue on the CPU time required and the efficiency of these methods. High order DC methods should be competitive in situations where using fully implicit methods is unavoidable.

\section{Conclusions}
\label{sec:conclusions}
We have presented a new approach of deferred correction methods for the numerical solution of general first order ordinary differential equations. Proofs for consistency, order of convergence and stability of the method are given, which rely on a recursive argument using a new deferred correction condition. The numerical experiments comply with the theory and show a high accuracy of the method and its satisfactory A-stable property and B-convergence. Globally, each DC scheme reaches its proper order of convergence and applies to any category of problem, providing accurate approximations for time steps not necessarily small. The accuracy of the DC schemes increases with the level of correction.



\bibliographystyle{spmpsci}      
\bibliography{references}   

\begin{thebibliography}{10}
\providecommand{\url}[1]{{#1}}
\providecommand{\urlprefix}{URL }
\expandafter\ifx\csname urlstyle\endcsname\relax
  \providecommand{\doi}[1]{DOI~\discretionary{}{}{}#1}\else
  \providecommand{\doi}{DOI~\discretionary{}{}{}\begingroup
  \urlstyle{rm}\Url}\fi

\bibitem{auzinger2016encyclopedia}
Auzinger, W.: Encyclopedia of {A}pplied and {C}omputational {M}athematics,
  chap. Defect Correction Methods, pp. 323--332.
\newblock Springer, Berlin, Heidelberg (2015)

\bibitem{IntegralDC2010}
Christlieb, A., Ong, B., Qiu, J.M.: Integral deferred correction methods
  constructed with high order {R}unge-{K}utta integrators.
\newblock Math. Comp. \textbf{79}, 761--783 (2010)

\bibitem{chung2010computational}
Chung, T.: Computational {F}luid {D}ynamics, 2nd edn.
\newblock Cambridge university press (2010)

\bibitem{dahlquist2008numerical}
Dahlquist, G., Bj\"{o}rck, A.k.: Numerical methods in scientific computing.
  {V}ol. {I}.
\newblock SIAM, Philadelphia, PA (2008)

\bibitem{MR0170477}
Dahlquist, G.G.: A special stability problem for linear multistep methods.
\newblock Nordisk Tidskr. Informationsbehandling (BIT) \textbf{3}, 27--43
  (1963)

\bibitem{daniel1967interated}
Daniel, J.W., Pereyra, V., Schumaker, L.L.: Iterated deferred corrections for
  initial value problems.
\newblock Acta Cient. Venezolana \textbf{19}, 128--135 (1968)

\bibitem{dutt2000spectral}
Dutt, A., Greengard, L., Rokhlin, V.: Spectral deferred correction methods for
  ordinary differential equations.
\newblock BIT \textbf{40}, 241--266 (2000)

\bibitem{enright1975comparing}
Enright, W.H., Hull, T., Lindberg, B.: Comparing numerical methods for stiff
  systems of {ODE}:s.
\newblock BIT \textbf{15}, 1--48 (1975)

\bibitem{frank1981concept}
Frank, R., Schneid, J., Ueberhuber, C.W.: The concept of {B}-convergence.
\newblock SIAM J. Numer. Anal. \textbf{18}, 753--780 (1981)

\bibitem{MR2058857}
Gustafsson, B., Kress, W.: Deferred correction methods for initial value
  problems.
\newblock BIT \textbf{41}, 986--995 (2001)

\bibitem{wanner1991solving}
Hairer, E., Wanner, G.: Solving ordinary differential equations. {II}. {S}tiff
  and differential-algebraic problems, vol.~14.
\newblock Springer-Verlag, Berlin (1991)

\bibitem{hansen2011order}
Hansen, A.C., Strain, J.: On the order of deferred correction.
\newblock Appl. {N}umer. {M}ath. \textbf{61}, 961--973 (2011)

\bibitem{hildebrand1974introduction}
Hildebrand, F.B.: Introduction to {N}umerical {A}nalysis.
\newblock McGraw-Hill Book Co., New York-D\"{u}sseldorf-Johannesburg (1974)

\bibitem{hull1972comparing}
Hull, T.E., Enright, W.H., Fellen, B.M., Sedgwick, A.E.: Comparing numerical
  methods for ordinary differential equations.
\newblock SIAM J. Numer. Anal. \textbf{9}, 603--637 (1972)

\bibitem{isaacson1966analysis}
Isaacson, E., Keller, H.B.: Analysis of numerical methods.
\newblock John Wiley \& Sons, Inc., New York-London-Sydney (1966)

\bibitem{karouma2015class}
Karouma, A.: A class of contractivity preserving {H}ermite-{B}irkhoff-{T}aylor
  high order time discretization methods.
\newblock Ph.D. thesis, Universit{\'e} d'Ottawa/University of Ottawa (2015)

\bibitem{koyaguerebo2020finite}
Koyaguerebo-Im{\'e}, S.C.E., Bourgault, Y.: Finite difference and numerical
  differentiation: General formulae from deferred corrections.
\newblock arXiv preprint arXiv:2005.11754  (2020)

\bibitem{koyaguerebo2020unconditionally}
Koyaguerebo-Im{\'e}, S.C.R., Bourgault, Y.: Arbitrary high-order
  unconditionally stable methods for reaction-diffusion equations via deferred
  correction: Case of the implicit midpoint rule.
\newblock arXiv:2006.02962v2.  (2020)

\bibitem{kraaijevanger1985b}
Kraaijevanger, J.: B-convergence of the implicit midpoint rule and the
  trapezoidal rule.
\newblock BIT \textbf{25}, 652--666 (1985)

\bibitem{kress2002deferred}
Kress, W., Gustafsson, B.: Deferred correction methods for initial boundary
  value problems.
\newblock J. Sci Comput. \textbf{17}(1-4), 241--251 (2002)

\bibitem{kushnir2012highly}
Kushnir, D., Rokhlin, V.: A highly accurate solver for stiff ordinary
  differential equations.
\newblock SIAM {J}. {S}ci. {C}omput. \textbf{34}, A1296--A1315 (2012)

\bibitem{quarteroni2010}
Quarteroni, A., Sacco, R., Saleri, F.: Numerical mathematics, vol.~37, second
  edn.
\newblock Springer-Verlag, Berlin (2007)

\bibitem{schild1990gaussian}
Schild, K.H.: Gaussian collocation via defect correction.
\newblock Numer. {M}ath. \textbf{58}, 369--386 (1990)

\bibitem{shampine1981evaluation}
Shampine, L.F.: Evaluation of a test set for stiff {ODE} solvers.
\newblock ACM Trans. Math. Software \textbf{7}, 409--420 (1981)

\bibitem{spijker1996stiffness}
Spijker, M.N.: Stiffness in numerical initial-value problems.
\newblock J. Comput. Appl. Math. \textbf{72}, 393--406 (1996)

\bibitem{stewart1990avoiding}
Stewart, K.: Avoiding stability-induced inefficiencies in {BDF} methods.
\newblock J. Comput. Appl. Math. \textbf{29}, 357--367 (1990)

\bibitem{tuenter2006frobenius}
Tuenter, H.J.: The {F}robenius problem, sums of powers of integers, and
  recurrences for the bernoulli numbers.
\newblock J. {N}umber {T}heory \textbf{117}, 376--386 (2006)

\end{thebibliography}


\end{document}